\documentclass[11pt]{amsart}

\usepackage[T1]{fontenc}
\usepackage{lmodern}
\usepackage{microtype}
\usepackage{amsmath,amssymb,mathtools}
\usepackage{enumitem}
\usepackage{xcolor}

\numberwithin{equation}{section}
\usepackage[colorlinks, citecolor=blue, linkcolor=red]{hyperref}

\newtheorem{theorem}{Theorem}[section]
\newtheorem{proposition}[theorem]{Proposition}
\newtheorem{lemma}[theorem]{Lemma}
\newtheorem{corollary}[theorem]{Corollary}
\theoremstyle{remark}
\newtheorem{remark}[theorem]{Remark}

\newcommand{\Ric}{\operatorname{Ric}}
\newcommand{\Riem}{\operatorname{Riem}}
\newcommand{\Hess}{\operatorname{Hess}}
\newcommand{\vol}{\operatorname{Vol}}
\newcommand{\supp}{\operatorname{supp}}
\newcommand{\T}{\mathbb T}
\newcommand{\R}{\mathbb R}
\newcommand{\cS}{\mathfrak S}
\newcommand{\dd}{\,\mathrm d}
\newcommand{\dv}{\dd V}

\newcommand{\set}[1]{{\left\{#1\right\}}}
\newcommand{\pa}[1]{{\left(#1\right)}}

\title[Complete critical metrics for the quadratic scalar curvature functional]
{Complete critical metrics for the quadratic scalar curvature functional: examples and rigidity}

\author{Giovanni Catino}
\address{Dipartimento di Matematica, Politecnico di Milano,
Piazza Leonardo da Vinci 32, 20133 Milano, Italy}
\email{giovanni.catino@polimi.it}

\author{Paolo Mastrolia}
\address{Dipartimento di Matematica, Universit\`a degli Studi di Milano,
Via Saldini 50, 20133 Milano, Italy}
\email{paolo.mastrolia@unimi.it}

\author{Dario D. Monticelli}
\address{Dipartimento di Matematica, Politecnico di Milano,
Piazza Leonardo da Vinci 32, 20133 Milano, Italy}
\email{dario.monticelli@polimi.it}


\hypersetup{
  pdftitle={Complete critical metrics for the quadratic scalar curvature functional: examples and rigidity},
  pdfauthor={Giovanni Catino, Paolo Mastrolia, Dario D. Monticelli},
  pdfsubject={Complete critical metrics for the quadratic scalar curvature functional},
  pdfkeywords={quadratic scalar curvature functional, critical metrics, rigidity, finite energy, scalar curvature}
}

\begin{document}

\begin{abstract}

We study complete critical metrics for the quadratic scalar curvature
functional
\[
g\mapsto\int_M S_g^2\,\dv_g,
\] with criticality understood
under compactly supported variations, focusing on rigidity phenomena and nontrivial examples. In dimension three, we prove that every complete
critical metric is scalar-flat, without any additional
assumption. In dimension four, we construct complete non-Einstein
critical metrics with bounded, nonconstant negative scalar curvature.
In every dimension $n\geq5$, we show that either $S_g\equiv 0$ or $S_g<0$ everywhere, and a new argument proves scalar-flatness whenever $S_g \in L^p$ for some $0<p\leq P_n$, where $P_n>2$ is an explicit threshold.
Moreover, we construct complete critical metrics with bounded, negative scalar curvature, decaying to zero at infinity.
These examples establish the optimality of our rigidity threshold $P_n$
in dimension six and in every dimension $n\geq10$. To the best of our knowledge, these are the first complete critical metrics for this functional that are neither scalar-flat nor Einstein.
\end{abstract}

\maketitle

\section{Introduction}

The aim of this paper is to continue the study of complete, non-compact critical metrics of the quadratic scalar curvature
functional
\[
        \cS^2(g)=\int_M S_g^2\,\dv_g.
\]

On a noncompact manifold, criticality is understood with respect to
compactly supported metric variations, without a volume constraint.
More precisely, given a smooth family of metrics $g_t$ with
$g_{t=0}=g$ and $g_t=g$ outside a fixed compact set, let $U\Subset M$
be a smooth relatively compact domain containing that set in its
interior. We require
\[
 \left.\frac{d}{dt}\right|_{t=0}\int_U S_{g_t}^2\,\dv_{g_t}=0;
\]
this condition is independent of the choice of $U$ and is equivalent
to the Euler--Lagrange equation below.

Throughout the paper manifolds are smooth, connected, and without
boundary, unless otherwise stated; statements for disconnected
manifolds apply to each connected component.

To fix the notation, let $M^{n}$, $n\geq 2$, be an $n$-dimensional smooth manifold without boundary. Given a Riemannian metric $g$ on $M^n$, we denote by $\Riem_g$, $W_g$, $\Ric_g$ and $S_g$ the Riemann curvature tensor, the Weyl tensor, the Ricci tensor and the scalar curvature, respectively.
The  \textit{quadratic curvature functionals} are
$$
\mathfrak{W}^2 = \int |W_g|^{2} dV_g\,, \qquad \mathfrak{r}^2 = \int |\Ric_g|^{2} dV_g\,, \qquad\, \mathfrak{S}^2 = \int S_g^{2} dV_g.
$$
For the corresponding pointwise quadratic curvature integrands, the
space is one-dimensional when $n=2$ (generated by $S_g^2$), and
two-dimensional when $n=3$ (generated by $S_g^2$ and $|\Ric_g|^2$).

From the standard decomposition of the Riemann tensor, for every $n\geq4$ one has
$$
\mathfrak{R}^2 = \int |\Riem_g|^{2} dV_g = \int \left(|W_g|^{2}+\frac{4}{n-2}|\Ric_g|^{2}-\frac{2}{(n-1)(n-2)}S_g^{2}\right) dV_g \,.
$$
 Quadratic curvature functionals have attracted a lot of attention from mathematicians and physicists in recent years (see e.g. \cite{berger3, Besse, gurvia3, CMM2023} and references therein). In \cite{CMM2023} we studied $\cS^2$ and the mixed quadratic functionals
\[
 \mathfrak F_t^2(g)=\int_M\bigl(|\Ric_g|^2+tS_g^2\bigr)\dv_g,
 \qquad t\in\R.
\]
In dimension three, every complete critical metric of $\mathfrak F_t^2$
is flat for $t>-1/3$, see \cite[Theorem~1.2]{CMM2023}. The homogeneous
examples of \cite{BGC2021} show that this threshold is sharp.

As far as the functional $\mathfrak{S}^2$ is concerned, in \cite{CMM2023} we proved the following for complete metrics:
 \begin{itemize}
   \item when $n=2$ all critical metrics are flat, and thus they are global minima of the functional;
   \item when $n=3$, all critical metrics are scalar-flat (and thus they are global minima of the functional) under the additional hypothesis that $S_g\in L^q(M^3)$ for some $q\in(1,\infty)$;
   \item the dimension $n=4$ is special, due to the differing behaviour of the quadratic functionals under constant rescaling of the metric when $n=4$ and $n\neq 4$; note also that, in this case, critical metrics turn out to have harmonic scalar curvature. Thus, if $M^4$ is compact, then it is scalar flat, or it has constant scalar curvature and it is Einstein; on the other hand, if $M^4$ is complete, non-compact and $S_g\in L^q(M^4)$ for some $q\in(1,\infty)$, then a classical result of Yau (see \cite{Yau1976}) implies that $M^4$ has constant scalar curvature, and hence it is either scalar-flat or a finite-volume Einstein metric;
   \item when $n\geq5$, there exists $q^*>2$ such that a critical metric of $\mathfrak{S}^2$ having scalar curvature $S_g$ which is \emph{bounded from below} and satisfying $S_g\in L^q(M^n)$ for $q\in(1,q^*)$ must be scalar-flat, and thus it is a global minimum of the functional.
 \end{itemize}
In \cite{CMM2026}, we removed the boundedness from below of the scalar curvature, assuming $n\geq 10$: the proof uses a conformal reduction to an auxiliary geometric system, in the same spirit as the works \cite{CMRgafa24, magliaro}. The case of the ``missing dimensions'' $5\leq n \leq 9$ has been subsequently settled by Zhang in \cite{Zhang2026}, thus completing the picture for finite energy critical metrics in dimensions greater than or equal to five.

We recall that the Euler--Lagrange equation for a critical metric $g$ of $\cS^2$ reads
\begin{equation}\label{eq:EL-full}
  2S_g\Ric_g-2\Hess(S_g)+2\Delta_gS_g\,g
       =\frac12S_g^2g.
\end{equation}
Equivalently,
\begin{align}
 S_g\Ric_g-\Hess(S_g)
    &=\frac{3}{4(n-1)}S_g^2g,\label{eq:EL-tensor}\\
 \Delta_gS_g
    &=\frac{n-4}{4(n-1)}S_g^2.\label{eq:EL-trace}
\end{align}
The classical
variational formulas may be found in
\cite[Proposition~4.66]{Besse}; for critical metrics in the complete setting
we refer to \cite{Catino2014,CMM2023,CMM2026}, and for the wider theory of
quadratic curvature energies to \cite{GurskyViaclovsky}.

The traced equation separates the compact and noncompact problems:
on a compact manifold with $n\neq4$, integration immediately gives
$S_g\equiv0$; in dimension four, $S_g$ is harmonic, hence constant
in the compact case; the tensor equation then gives the alternative
between scalar-flatness and the Einstein condition. On a complete
noncompact manifold, these integrations require control at infinity.
The results below address both the necessity of such control and the
range of integrability assumptions that yields rigidity.

In \cite{Catino2014}, complete critical metrics with nonnegative
scalar curvature were classified, and
the necessity of the sign assumption was left open. Our first two
results answer this question in opposite ways: in dimension three
the sign assumption can be removed entirely, whereas in dimension
four there exist complete critical metrics with bounded, nonconstant
negative scalar curvature.

In dimension three, the following theorem also removes the
integrability assumption from
\cite[Theorem~1.3 and Remark~1.4]{CMM2023}:

\begin{theorem} \label{thm:three-dimensional-rigidityINTRO}
Let $(M^3,g)$ be a three-dimensional complete critical metric of
$\cS^2$. Then $(M^3,g)$ is scalar-flat, and hence a global
minimum of the functional.
\end{theorem}

In dimension four, completeness alone no longer gives the same rigidity: indeed, the following examples have bounded, nonconstant negative scalar
curvature and are neither scalar-flat nor Einstein. They also show
why an additional hypothesis remains necessary in the $L^p$
classification of Proposition~\ref{prop:four-lp-classification} below:
their harmonic scalar curvature belongs to no finite $L^p$ space.

\begin{theorem}\label{thm:four-dimensional-existenceINTRO}
There exists a smooth, complete, Riemannian metric $g_4$ on
$M^4=\R\times N$, where $N$ is a closed hyperbolic three-manifold,
which is critical for $\cS^2$ and has the following properties:
\begin{enumerate}[label=\textup{(\roman*)}]
\item $g_4$ is neither scalar-flat nor Einstein;
\item $S_{g_4}$ is negative, bounded, and nonconstant on $M^4$;
\item for every $q>0$,
\[
 S_{g_4}\not\in L^q(M^4,g_4).
\]
In particular, $g_4$ has infinite energy.
\end{enumerate}
\end{theorem}

In the second part of the paper we treat all dimensions $n\geq5$, by a new, unified argument
in the original metric (without any conformal change of $g$). Besides recovering finite-energy rigidity,
the theorem establishes rigidity in the range $0<p\leq n/2$ for $5\leq n\leq9$
and $0<p\leq\frac{7n-10}{2(n-4)}$ for $n\geq10$, thus including the endpoint in both cases.

\begin{theorem}\label{thm:rigidity-intro}
Let $(M^n,g)$, $n\geq5$, be a complete critical metric of $\cS^2$.
Then either $S_g\equiv0$ or $S_g<0$ everywhere.
Moreover, if $S_g\in L^p(M,g)$ for some $0<p\leq P_n$, with
\begin{equation}\label{eq:p-range-intro}
 P_n=\begin{cases}
 n/2,&5\leq n\leq9,\\
 \frac{7n-10}{2(n-4)},&n\geq10,
 \end{cases}
\end{equation}
then $S_g\equiv0$. In particular, for $n\geq5$, every complete
critical metric of $\cS^2$ with finite energy is scalar-flat.
\end{theorem}

The sign alternative in Theorem~\ref{thm:rigidity-intro} is proved
in Lemma~\ref{lem:unconditional-sign}, independently of any
integrability assumption. It excludes sign-changing scalar curvature
for every complete critical metric in dimension $n\geq5$.
In the negative case, set $u=-S_g>0$ and $w=u^{-1/2}$. A weighted Bochner argument with
respect to the finite measure $u^p\,\dv_g$ first bounds
$|\nabla w|^2$ by a constant. A second cutoff test yields the dimensional gradient bound. The scalar equation and a cutoff argument prove rigidity
below the upper endpoint; at equality, they force $|\nabla w|$ to be a
positive constant, which contradicts completeness and positivity of $w$.

Finally, in the third part of the paper we construct a family of explicit complete critical metrics on
$\R^2\times\T^{n-2}$, where $\T^{n-2}$ is a flat torus: their scalar curvature is bounded and tends to zero at infinity,
but is not identically zero; moreover, the exact curvature integrals also show that
the rigidity endpoint is optimal in every dimension $n\geq10$.

\begin{theorem}\label{thm:examples-intro}
For every $n\geq5$, there exists a smooth complete Riemannian
metric \(g_n\) on
\[
              M_n=\R^2\times\T^{n-2}
\]
which is critical for $\cS^2$ and has the following properties:
\begin{enumerate}[label=\textup{(\roman*)}]
\item \(g_n\) is not scalar-flat;
\item \(-1\leq S_{g_n}<0\) on \(M_n\), and $S_{g_n}$ tends to zero at infinity;
\item for every \(q>0\),
\[
 S_{g_n}\in L^q(M_n,g_n)
 \quad\Longleftrightarrow\quad
 q>q_n=\frac{7n-10}{2(n-4)}.
\]
\end{enumerate}
In particular, $g_n$ has infinite energy.
\end{theorem}

The previous theorem thus implies that completeness, bounded scalar curvature, and even pointwise decay to zero
do not force scalar-flatness. Moreover,
\[
            q_n-2=\frac{3(n+2)}{2(n-4)}>0,
\]
so none of these metrics has finite energy. When
\(n\geq10\), this highlights sharpness of Theorem~\ref{thm:rigidity-intro}.
A second construction in dimension $n=6$ shows that the rigidity threshold \(P_6=3\) is optimal.

\begin{theorem}\label{thm:n6-sharp-intro}
Let $(N^5,\bar h)$ be closed and Einstein, with
$\Ric_{\bar h}=-\bar h/10$. On $M = \R\times N$, set
\[
 \widehat g_6=\cosh(t/\sqrt{10})\,(dt^2+\bar h).
\]
Then $\widehat g_6$ is a smooth complete non-Einstein critical metric, with
\[
 S_{\widehat g_6}=-\operatorname{sech}(t/\sqrt{10}),
 \qquad
 S_{\widehat g_6}\in L^p(M,\widehat g_6)\quad\Longleftrightarrow\quad p>3.
\]
\end{theorem}

\begin{remark}

Theorem~\ref{thm:n6-sharp-intro} shows that the endpoint $P_6=3$ in
Theorem~\ref{thm:rigidity-intro} is optimal. The change of variables
establishing the exact integrability range is given in the proof of
Proposition~\ref{prop:n6-sharp-example}.

\end{remark}

In dimension six and in dimension $n\geq 10$, the examples also show that the constant in the gradient estimate of Proposition~\ref{prop:gradient-estimate} is
optimal, see Remark~\ref{rem:sharp-gradient-examples}. In dimensions $5,7,8,9$, it remains open whether the rigidity range
in Theorem~\ref{thm:rigidity-intro} can be enlarged.

The paper is organized as follows. Sections~\ref{sec:three-dimensional-case}
and~\ref{sec:four-dimensional-case} concern dimensions three and four.
Section~\ref{sec:direct-rigidity} proves rigidity for
$n\geq5$. In Section~\ref{sec:direct-examples} we construct the family of Theorem \ref{thm:examples-intro} and the six-dimensional example of Theorem \ref{thm:n6-sharp-intro}, and establish the sharpness claims.

\medskip

\emph{Acknowledgements.} The first and second authors are members of the GNSAGA, Gruppo Nazionale per le Strutture Algebriche, Geometriche e le
loro Applicazioni of INdAM. The third author is a member of GNAMPA,
Gruppo Nazionale per l’Analisi Matematica, la Probabilità e le loro Applicazioni of INdAM.

\section{Dimension three: rigidity}
\label{sec:three-dimensional-case}
In this section we prove Theorem~\ref{thm:three-dimensional-rigidityINTRO}.
Recall that a \emph{nodal domain} of a continuous function $f$ is a
connected component of the open set $\{f\neq0\}$. On each nodal domain,
$f$ has a fixed strict sign. Accordingly, a \emph{negative} (respectively,
\emph{positive}) nodal domain of $f$ is a connected component of
$\{f<0\}$ (respectively, $\{f>0\}$).

The critical equations in dimension three, dropping the subscript $g$ for the sake of readability, are
\begin{align}
 S\Ric-\Hess(S)&=\frac38S^2g,
 \label{eq:three-critical-tensor}\\
 \Delta S&=-\frac18S^2.
 \label{eq:three-critical-trace}
\end{align}
The case $S\geq0$ is covered in \cite[Theorem~1.2]{Catino2014}; it therefore suffices to exclude a connected component
$\Omega\subset\{S<0\}$. On such a component, $u=-S>0$ satisfies
\begin{align}
 \Hess(u)&=u\Ric+\frac38u^2g,
 \label{eq:three-u-hessian}\\
 \Delta u&=\frac18u^2.
 \label{eq:three-u-laplacian}
\end{align}

The proof has three steps. First, the boundary of a negative nodal
domain is regular and totally geodesic, and its intrinsic completion
has minimizing segments that stay in the interior. Second, the tensor
equation yields a local weighted distance comparison and cutoffs in
the barrier sense. Finally, point selection and the maximum principle
applied to a smooth lower support give a contradiction.

\subsection*{Step 1.} Since $\Omega$ may not be complete, we use the length distance of the
restricted metric, namely
\[
 d_\Omega(x,y)=\inf_c\int_0^1|\dot c(t)|_g\,dt,
\]
where the infimum is over piecewise smooth curves $c:[0,1]\to\Omega$
joining $x$ to $y$. The curves are required to stay in $\Omega$, so
this distance may exceed the ambient distance $d_M$.
By the \emph{intrinsic completion} we mean the ordinary metric
completion $X$ of $(\Omega,d_\Omega)$: its points are equivalence
classes of $d_\Omega$-Cauchy sequences, and
\[
 d_X([(x_j)],[(y_j)])=\lim_{j\to\infty}d_\Omega(x_j,y_j).
\]
We identify $\Omega$ with its isometric image in $X$. A \emph{length space}
is a metric space whose distance is the infimum of curve lengths;
it is called \emph{geodesic} if this infimum is attained for every pair of
points, and proper if every closed bounded ball is compact.
For these notions see e.g. \cite[Sections~1.5 and~2.1--2.5]{BuragoBuragoIvanov}.
The next lemma proves the required properties of $X$ and identifies
it with the ambient closure of $\Omega$. All subsequent differential inequalities are used only in the smooth
interior $\Omega$.

\begin{lemma}\label{lem:nodal-set}
Suppose that $S\not\equiv0$. Then $0$ is a regular value of $S$,
and $\partial\Omega$ is totally geodesic.
The intrinsic completion $(X,d_X)$ of $(\Omega,g)$ is naturally
homeomorphic to $\overline\Omega\subset M$ and is proper and geodesic.
Every minimizing segment in $X$ whose endpoints belong to $\Omega$
is contained in $\Omega$.
\end{lemma}

\begin{proof}
Suppose that $S(x)=0$ and $\nabla S(x)=0$. Along every unit-speed
geodesic $\gamma$ originating from $x$, the function $v=S\circ\gamma$ satisfies
\[
 v''=\Ric(\dot\gamma,\dot\gamma)v-\frac38v^2,
 \qquad v(0)=v'(0)=0.
\]
Standard uniqueness for this ordinary differential equation gives $S=0$ in
a normal neighborhood of $x$; thus the set $\{S=0, \nabla S=0\}$ is both
open and closed. Connectedness and $S\not\equiv0$ imply that this
set is empty. At a zero of $S$, the tensor equation also gives
$\Hess(S)=0$. The second fundamental form of the regular level
$\{S=0\}$ is, up to the choice of the normal,
$|\nabla S|^{-1}\Hess(S)$ restricted to its tangent bundle;
it therefore vanishes.

Let $d_\Omega$ denote the intrinsic distance. Since
$d_M\leq d_\Omega$, every intrinsic Cauchy sequence has an ambient
limit in $\overline\Omega$, by completeness of $M$.
At an interior limit point, a small coordinate ball identifies
intrinsic and ambient convergence. At a boundary point, regularity
of the zero level gives a coordinate half-ball on the negative side.
This negative side belongs to $\Omega$, since it is connected and
meets $\Omega$. Local upper and lower bounds for the metric
coefficients show that two sequences approaching the same boundary
point have intrinsic distance tending to zero. They also show that
every boundary point is the limit of an intrinsic Cauchy sequence.
These observations identify $X$ with $\overline\Omega$ and establish
equality of the two topologies.
Under this identification, $d_X$ is also the length distance of
the metric on the manifold with boundary $\overline\Omega$.
Indeed, curves in $\Omega$ are admissible in $\overline\Omega$.
Conversely, a piecewise smooth curve in $\overline\Omega$ can be
approximated, with converging lengths, by curves in $\Omega$:
in boundary half-charts move its boundary portions a distance
$\varepsilon$ into the interior, using a smooth inward-pointing
vector field on a neighborhood of the compact image of the curve.
The flow tends to the identity in $C^1$ as $\varepsilon\to0$;
short interior arcs restore any prescribed interior endpoints.
Passing to the completion gives the same equality when an endpoint
lies on the boundary. This is an equality with the intrinsic length
distance of $\overline\Omega$, not with the ambient distance $d_M$.
The half-ball description proves local compactness of $X$.
Moreover, the completion of a length space is a length space.
The Hopf--Rinow theorem for length spaces now implies that $X$ is
proper, meaning that its closed bounded balls are compact, and that
any two points can be joined by a minimizing segment (see
\cite[Proposition~2.5.22 and Theorem~2.5.23]{BuragoBuragoIvanov}).

It remains to show that minimizing segments with interior endpoints
avoid the boundary. If $\partial\Omega=\varnothing$, the assertion
is immediate. Otherwise form the double $\widehat X$ by gluing two
copies of $X$ along their common boundary. In Fermi coordinates,
total geodesicity gives
\[
 g=dt^2+g_t,\qquad g_t=g_0+t^2A(t,z),\qquad t\geq0,
\]
with $A$ smooth up to $t=0$.
The reflected metric is
\[
 \widehat g=dt^2+g_0+t^2A(|t|,z).
\]
It is $C^2$ across $t=0$: since the second fundamental  form vanishes, the first normal derivative vanishes as well;
the second normal derivative tends to $2A(0,z)$ from both sides,
and mixed tangential-normal derivatives tend to zero. Tangential
derivatives have the same limits as well. Its Christoffel symbols
are therefore $C^1$, hence locally Lipschitz, and the geodesic
initial value problem has uniqueness. No smoothness of the doubled
metric beyond $C^2$ is needed.

Folding a curve in $\widehat X$ into one copy of $X$ does not increase
its length. Inclusion and folding therefore give equality of the
distances in $X$ and $\widehat X$ for points of the same copy.

An $X$-minimizing segment is consequently a length-minimizing
curve for the $C^2$ metric on the double. After constant-speed
parametrization it satisfies the geodesic equation; this follows
from the local first variation of the energy and requires only
the stated $C^2$ regularity.

If it touched $t=0$ while remaining in $t\geq0$, its normal velocity
would vanish at the contact point. Reflection would produce a
geodesic with the same position and velocity. Uniqueness would make
the whole segment fixed by reflection, hence contained in the
boundary, contradicting the interior endpoints.
\end{proof}

The boundary-avoidance property will also be used in the standard
\emph{Calabi trick}: a minimizing segment in $X$ with interior endpoints
can be shifted slightly at its starting point to produce a smooth
upper support for the intrinsic distance at its endpoint.

\subsection*{Step 2.} Extend $u$ continuously to $X$ by zero on $\partial\Omega$, and set
on $\Omega$
\[
 d\mu=u\,\dv,\qquad f=-\log u,\qquad
 L=\Delta-\langle\nabla f,\nabla\cdot\rangle
   =u^{-1}\operatorname{div}(u\nabla\cdot).
\]

We use the convention
\[
 \Ric_f^m=\Ric+\Hess(f)-\frac1m\,df\otimes df\qquad(m>0),
\]
where $m$ is the additional dimension parameter: the effective
dimension is $n+m$. In the present three-dimensional argument
$m=1$, so the effective dimension is four. This is the convention
used in \cite[Appendix~A]{WeiWylie}.
Since
\[
 \Hess(f)=-\frac{\Hess(u)}u+\frac{du\otimes du}{u^2},\qquad
 df\otimes df=\frac{du\otimes du}{u^2},
\]
the critical equations imply
\begin{equation}\label{eq:three-bakry-emery}
 Lu=\frac18u^2+\frac{|\nabla u|^2}{u}\geq\frac18u^2,
 \qquad
 \Ric_f^1=\Ric+\Hess(f)-df\otimes df=-\frac38u\,g.
\end{equation}

The weight can vanish at $\partial\Omega$, but only interior pointwise calculations
will be used.

\begin{lemma}
\label{lem:three-weighted-package}
There is a universal constant $C$ with the following property.
Suppose that $q\in\Omega$, $r>0$, and $u\leq\Lambda$ on
$B_{8r}^X(q)$. Then there is a Lipschitz function
$\eta:X\to[0,1]$, equal to one on $B_r^X(q)$ and supported
compactly in $B_{2r}^X(q)$, with the following property:
for every $x\in\Omega$ such that $\eta(x)>0$, there is a smooth
function $\eta_x$ in a neighborhood of $x$ with
\[
 \eta_x\leq\eta,\qquad \eta_x(x)=\eta(x)>0,
\]
and
\begin{equation}\label{eq:three-cutoff-laplacian}
 \frac{|\nabla\eta_x(x)|^2}{\eta(x)}\leq\frac{C}{r^2},
 \qquad L\eta_x(x)\geq
       -C\left(\frac1{r^2}+\frac{\sqrt\Lambda}{r}\right).
\end{equation}
\end{lemma}

\begin{proof}
Write $\rho=d_X(q,\cdot)$ and $\kappa=\sqrt{\Lambda/8}$.
Every minimizing segment from $q$ to an interior point of
$B_{3r}^X(q)$ stays in the smooth interior by
Lemma~\ref{lem:nodal-set}, and lies in $B_{3r}^X(q)$.
Along such a segment,
$\Ric_f^1=-3u\,g/8\geq-3\kappa^2g$.
At a point where $\rho$ is smooth, put $m_\rho=L\rho$.
Bochner's formula and the radial Hessian inequality give
\[
 m_\rho'=-|\Hess\rho|^2-\Ric_f^1(\dot\gamma,\dot\gamma)
            -\bigl((f\circ\gamma)'\bigr)^2
 \leq-\frac{m_\rho^2}{3}+3\kappa^2.
\]
Here the inequality uses
$|\Hess\rho|^2\geq(\Delta\rho)^2/2$ and
$(\Delta\rho)^2/2+(f')^2\geq(L\rho)^2/3$.
Since $m_\rho(t)=2/t+O(1)$ near the center, Riccati comparison
yields the following local form of
\cite[Theorem~A.1]{WeiWylie}:
\begin{equation}\label{eq:three-local-distance-comparison}
 L\rho\leq3\kappa\coth(\kappa\rho)
           \leq\frac3\rho+3\kappa
\end{equation}
up to the first cut point, with the interpretation $3/\rho$
when $\kappa=0$.
We now use this estimate at an interior point
$x\in B_{2r}^X(q)\setminus\{q\}$, which is the only region
needed for the cutoff.
Let $\gamma:[0,\ell]\to X$ be an arclength minimizing segment
from $q$ to $x$. It is contained in $\Omega$.
For $0<\delta<\ell$, set $q_\delta=\gamma(\delta)$ and
\[
 \rho_\delta(y)=\delta+d_X(q_\delta,y).
\]
The triangle inequality gives $\rho_\delta\geq\rho$ with
$\rho_\delta(x)=\rho(x)=\ell$.
Moreover $d_X(q_\delta,\cdot)$ is smooth near $x$.
Indeed, a second minimizing segment from $q_\delta$ to $x$
would concatenate with $\gamma|_{[0,\delta]}$ to give a
minimizing broken curve from $q$ to $x$, and therefore must
coincide with $\gamma$.
Nor can $x$ be conjugate to $q_\delta$ along $\gamma$.
To see this, suppose that a nonzero Jacobi field $J$ vanishes
at $\delta$ and $\ell$, and extend it by zero on $[0,\delta]$.
This is a continuous piecewise smooth field with zero endpoints.
The index form $I$ of $\gamma$ is nonnegative, because $\gamma$
minimizes and its image is a compact subset of $\Omega$.
The Jacobi equation gives $I(J,J)=0$. Applying nonnegativity to
$J+\tau V$, for every real $\tau$ and every smooth endpoint-zero
field $V$, then gives $I(J,V)=0$. Integration by parts on
$[\delta,\ell]$ yields
\[
 I(J,V)=-\langle D_tJ(\delta^+),V(\delta)\rangle.
\]
Since $V(\delta)$ is arbitrary, $D_tJ(\delta^+)=0$.
Together with $J(\delta)=0$, uniqueness for the Jacobi equation
would force $J=0$, a contradiction.
\par
Finally, choose a compact neighborhood $K\subset\Omega$ of $x$.
All minimizing segments from $q_\delta$ to points of $K$ have
uniformly bounded lengths and lie in a fixed compact ball of $X$.
If they were not uniformly separated from $X\setminus\Omega$,
a sequence of such segments, parametrized on $[0,1]$ with constant
speed, would converge uniformly to a minimizing segment from
$q_\delta$ to some point of $K$ that touches the boundary.
Here one uses continuity of the distance and lower semicontinuity
of length; the limit has the distance between its endpoints as
its length. Lemma~\ref{lem:nodal-set} excludes this possibility.
Thus these segments lie in a common compact subset of $\Omega$.

The exponential map at $q_\delta$ is defined on a neighborhood
of the initial vector of $\gamma|_{[\delta,\ell]}$, and its
differential is invertible there by the absence of conjugacy.
Every minimizing segment to an endpoint tending to $x$ has initial
vector tending to this vector: otherwise the preceding compactness
and the geodesic equation would give a different minimizing segment
to $x$. The local inverse of the exponential map therefore describes
all nearby minimizers, and their lengths are a smooth function of
the endpoint. This proves that $d_X(q_\delta,\cdot)$ is smooth
near $x$, despite the possible incompleteness of $\Omega$.
The segment from $q_\delta$ to $x$ is a subsegment of $\gamma$,
hence lies in $B_{2r}^X(q)$, where $u\leq\Lambda$.
The same Riccati comparison from $q_\delta$ to $x$ therefore gives

\begin{equation}\label{eq:calabi-weighted-shift}
 L\rho_\delta(x)\leq
     \frac{3}{\ell-\delta}+3\kappa.
\end{equation}

Choose smooth nonincreasing $\zeta:[0,\infty)\to[0,1]$,
equal to one on $[0,1]$, compactly supported in $[0,2)$,
and satisfying $(\zeta')^2\leq C\zeta$ and
$|\zeta''|\leq C$. Set $\eta=\zeta(\rho/r)$.
This function is Lipschitz and has compact support because $X$
is proper.
When $\rho(x)<r$, one can take $\eta_x=1$.
For $r\leq\rho(x)<2r$ and $\eta(x)>0$, take
\[
 \eta_x(y)=\zeta\bigl(\rho_\delta(y)/r\bigr)
\]
near $x$, with $0<\delta<r/2$.
Since $\rho_\delta\geq\rho$ and $\zeta'\leq0$,
this is a smooth lower support for $\eta$ at $x$.
At $x$, $|\nabla\rho_\delta|=1$ and the chain rule, together
with \eqref{eq:calabi-weighted-shift} and
$\rho(x)-\delta\geq r/2$, yields
\[
 \frac{|\nabla\eta_x|^2}{\eta(x)}\leq\frac C{r^2},
 \qquad
 L\eta_x(x)\geq-C\left(\frac1{r^2}+\frac\kappa r\right).
\]
Increasing the universal constant proves the assertion.
\end{proof}

\subsection*{Step 3.} Assume that a negative nodal domain $\Omega$ exists, and retain the
notation above. Fix $x_0\in\Omega$ and $A>1$. We first choose a
point $q\in\Omega$ at whose natural scale $u$ is controlled:
with $U=u(q)>0$, we require
\begin{equation}\label{eq:three-pointwise-doubling}
 U\geq u(x_0),\qquad
 u\leq4U\ \text{on }B_{8A/\sqrt U}^X(q).
\end{equation}

The following point-selection argument is the doubling construction
of \cite[Lemma~5.1]{PolacikQuittnerSouplet}, adapted here to the
continuous nonnegative function $\sqrt u$ on the complete space $X$.
We include the argument so that its behavior at the nodal boundary
is explicit.
Start at $q_0=x_0$. If the desired local bound fails at $q_j$, choose

$q_{j+1}$ within distance $8A/\sqrt{u(q_j)}$ where
$\sqrt{u(q_{j+1})}>2\sqrt{u(q_j)}$. An infinite sequence would satisfy
\[
 \sqrt{u(q_j)}\geq2^j\sqrt{u(x_0)},\qquad
 \sum_{j=0}^{\infty}d_X(q_j,q_{j+1})
                 \leq\frac{16A}{\sqrt{u(x_0)}}.
\]
It would therefore converge in the complete space $X$, contradicting
continuity of the extension of $u$ because $u(q_j)\to\infty$.
The selection must stop, proving \eqref{eq:three-pointwise-doubling}.

Apply Lemma~\ref{lem:three-weighted-package} with
$r=A/\sqrt U$ and $\Lambda=4U$. It supplies a cutoff $\eta$
whose smooth lower supports satisfy
\[
 \frac{|\nabla\eta_x(x)|^2}{\eta(x)}\leq\frac{CU}{A^2},
 \qquad L\eta_x(x)\geq-CU(A^{-1}+A^{-2}).
\]
The function $F=\eta u$, extended by zero to $X\setminus\Omega$,
is continuous and compactly supported in $X$, since $X$ is proper
and $u=0$ on $\partial\Omega$. Therefore it attains its
positive maximum at some $x_*\in\Omega$, with
$F(x_*)\geq F(q)=U$ and $\eta(x_*)>0$.
Let $\eta_*$ be the smooth lower support given by the lemma at
$x_*$. Since $u>0$ nearby, $F_*=u\eta_*$ has a local maximum
at $x_*$, so
\[
 \nabla F_*(x_*)=0,\qquad LF_*(x_*)\leq0.
\]
The first identity gives
$\nabla u=-u\nabla\eta_*/\eta_*$ at $x_*$. Hence
\begin{align*}
 0&\geq LF_*(x_*)
 =\eta_*Lu+uL\eta_*+2\langle\nabla u,\nabla\eta_*\rangle\\
 &\geq u\left(\frac{F(x_*)}{8}
           -CU(A^{-1}+A^{-2})
           -\frac{2CU}{A^2}\right).
\end{align*}
Choosing $A$ universally large makes the parentheses strictly
positive, because $F(x_*)\geq U$. This contradiction excludes every
negative nodal domain.
Therefore $S\geq0$, and \cite[Theorem~1.2]{Catino2014} gives
$S\equiv0$
and this concludes the proof of Theorem~\ref{thm:three-dimensional-rigidityINTRO}.

\section{Dimension four: complete non-Einstein examples}
\label{sec:four-dimensional-case}

In dimension four, the critical equations become
\begin{equation}
 S\Ric-\Hess(S)=\frac14S^2g,\qquad \Delta S=0.
\label{eq:four-critical-system}
\end{equation}
As already recalled in the introduction, harmonicity settles the $L^p$ case. The construction below shows
that bounded scalar curvature alone does not imply rigidity.

\begin{proposition}
\label{prop:four-lp-classification}
Let $(M^4,g)$ be complete, connected, without boundary, and critical
for $\cS^2$. If $M$ is compact, or if $S\in L^p(M,g)$ for some
$1<p<\infty$, then $S$ is constant and $g$ is scalar-flat or Einstein.
In the noncompact case, a nonzero constant scalar curvature under
this $L^p$ assumption is necessarily negative and $\vol_g(M)<\infty$.
\end{proposition}

\begin{proof}
On a compact manifold, harmonicity forces $S$ to be constant; in the complete
case, the same result follows by the classical $L^p$ Liouville argument of Yau \cite{Yau1976}.

If $S=c\neq0$, \eqref{eq:four-critical-system} gives
$\Ric=(c/4)g$. In the noncompact case the $L^p$ assumption implies
$\vol_g(M)<\infty$, and the
Bonnet--Myers theorem forces $c<0$.
\end{proof}

The Einstein alternative is necessary: a complete finite-volume
hyperbolic four-manifold has $S=-12\in L^p$ for every $p>0$.
Also, a complete noncompact critical metric with $S\geq0$ is
scalar-flat by \cite[Theorem~1.3]{Catino2014} and the Bonnet--Myers theorem.
The following family has bounded, strictly negative, nonconstant
scalar curvature with infinite energy.

\begin{proposition}[Complete non-Einstein examples]
\label{thm:four-complete-examples}
Let $(N^3,\bar h)$ be a closed connected hyperbolic manifold with
$\Ric_{\bar h}=-2\bar h$. For $0<Q<16/\sqrt3$, set
\[
 V_Q(r)=-1+\frac{r^2}{12}+\frac{Q^2}{4r^4},\qquad r>0,
\]
and let $b_0$ be the larger positive zero of $V_Q$. The solution of
\begin{equation}
 b''=\frac b{12}-\frac{Q^2}{2b^5},\qquad
 b(0)=b_0,\qquad b'(0)=0
\label{eq:four-b-ode}
\end{equation}
is positive, even, and defined on $\R$. Define
\[
 w(s)=Q\int_0^s b(t)^{-3}\,dt,\qquad
 W_Q=Q\int_0^\infty b(s)^{-3}\,ds.
\]
Then $0<W_Q<\infty$, and
\begin{equation}
 g_Q=e^{-w(s)}\bigl(ds^2+b(s)^2{\bar h}\bigr)
\label{eq:four-example-metric}
\end{equation}
is a smooth complete critical metric of $\cS^2$ on $\R\times N$.
It is not Einstein, and its scalar curvature satisfies
\[
 S_{g_Q}=-e^{w(s)},\qquad
 \lim_{s\to\pm\infty}S_{g_Q}(s)=-e^{\pm W_Q}.
\]
In particular, $S_{g_Q}$ is bounded and belongs to no
$L^p(M,g_Q)$ with $0<p<\infty$.
\end{proposition}

\begin{proof}
For $h=ds^2+b^2\bar h$ and $w=w(s)$, consider the auxiliary system
\begin{equation}
 \Ric_h=-\frac14h+\frac32\,dw\otimes dw,\qquad \Delta_hw=0.
 \label{eq:four-auxiliary-system}
\end{equation}
We will verify below that this system makes $e^{-w}h$ a critical
metric. The warped-product formulas are
\[
 \Ric_h(\partial_s,\partial_s)=-3\frac{b''}{b},\qquad
 \left.\Ric_h\right|_{TN}
 =\frac{-2-bb''-2(b')^2}{b^2}\left.h\right|_{TN},\qquad
 \Delta_hw=\frac{(b^3w')'}{b^3}.
\]
Thus, the harmonicity of $w$ gives $b^3w'=Q$, with $Q$ constant. The radial
Ricci equation becomes
\[
 -3\frac{b''}{b}=-\frac14+\frac{3Q^2}{2b^6},
 \qquad\text{hence}\qquad
 b''=\frac b{12}-\frac{Q^2}{2b^5}.
\]
The tangential equation is $-2-bb''-2(b')^2=-b^2/4$.
Substituting the expression for $b''$ yields
\[
 (b')^2=-1+\frac{b^2}{12}+\frac{Q^2}{4b^4}=V_Q(b).
\]
This derives both the ODE and its first integral from the auxiliary
system. Conversely, these two equations together with $w'=Qb^{-3}$
imply \eqref{eq:four-auxiliary-system}.

The function $V_Q$ tends to $+\infty$ as $r\to0^+$ and as $r\to\infty$, and has
a unique critical point, at which
\[
 r_*=(6Q^2)^{1/6},\qquad
 V_Q(r_*)=-1+\frac{(6Q^2)^{1/3}}8<0.
\]
Thus, it has two simple positive zeros. Since
$V_Q(\sqrt8)<0<V_Q(\sqrt{12})$, the larger zero satisfies
$\sqrt8<b_0<\sqrt{12}$ and $V_Q'(b_0)>0$.
The initial value problem and its first integral give
\begin{equation}
 (b')^2=V_Q(b),\qquad b''(0)=\tfrac12V_Q'(b_0)>0.
\label{eq:four-b-first-integral}
\end{equation}

Uniqueness makes $b$ even. In its positive-time interval of
existence, $b'>0$ and the first integral give
$b'=\sqrt{V_Q(b)}$. Separating variables yields
\[
 s=\int_{b_0}^{b(s)}\frac{dr}{\sqrt{V_Q(r)}}.
\]
This integral is finite at its lower endpoint because $b_0$ is a
simple zero and $V_Q'(b_0)>0$. At infinity its integrand is
asymptotic to $\sqrt{12}/r$, so
\[
 \int_{b_0}^{\infty}\frac{dr}{\sqrt{V_Q(r)}}=\infty.
\]
Thus $b$ cannot tend to infinity at a finite value of $s$.
It also cannot approach a finite limit larger than $b_0$, since
$b'=\sqrt{V_Q(b)}>0$ there. The ODE continuation theorem now gives
existence for all $s\geq0$ and $b(s)\to\infty$; evenness gives the
same conclusion at the other end.

The functions just constructed define $h,w$ on  $M=\R\times N$ and they satisfy
\eqref{eq:four-auxiliary-system} on  $M$.
Let $u=e^w$, so that $h=ug_Q$.
Since $h=e^wg_Q$, the four-dimensional conformal formulas give
\[
 \Delta_hw=e^{-w}\bigl(\Delta_{g_Q}w+|\nabla w|_{g_Q}^2\bigr).
\]
Thus $\Delta_hw=0$ implies both
$\Delta_{g_Q}w+|\nabla w|_{g_Q}^2=0$ and
$\Delta_{g_Q}u=u(\Delta_{g_Q}w+|\nabla w|_{g_Q}^2)=0$.
The  transformation law for the Ricci tensor  is (see e.g. \cite{CMBook})
\[
 \Ric_h=\Ric_{g_Q}-\Hess_{g_Q}(w)
 +\tfrac12dw\otimes dw
 -\tfrac12\bigl(\Delta_{g_Q}w+|\nabla w|_{g_Q}^2\bigr)g_Q.
\]
The last term vanishes by the identity just proved.

Comparing with \eqref{eq:four-auxiliary-system}, we obtain
\[
 \Hess_{g_Q}(u)=u\Ric_{g_Q}+\frac14u^2g_Q.
\]
Its trace gives $S_{g_Q}=-u$, and consequently
\eqref{eq:four-critical-system} holds.

Finally,
\[
 \int_0^\infty b(s)^{-3}\,ds
 =\int_{b_0}^\infty\frac{dr}{r^3\sqrt{V_Q(r)}}<\infty.
\]
Integrability follows from the simple zero at $b_0$ and decay
of order $r^{-4}$ at infinity. Hence $w$ is odd and strictly
increasing, with finite limits $\pm W_Q$, and bounded.
Since $b\geq b_0>0$ and $N$ is closed, $h$ is complete, and
\[
 e^{-W_Q}h\leq g_Q\leq e^{W_Q}h
\]
proves completeness of $g_Q$. The function $S_{g_Q}=-e^w$ is
nonconstant, so $g_Q$ is not Einstein. For every $p>0$,
\[
 \int_M|S_{g_Q}|^p\,\dv_{g_Q}
 =\vol_{\bar h}(N)\int_\R e^{(p-2)w(s)}b(s)^3\,ds=\infty,
\]
because $w$ is bounded and $b\geq b_0>0$.
\end{proof}

\section{Higher-dimensional rigidity}\label{sec:direct-rigidity}

In this section we first settle the sign of the scalar curvature without any
integrability assumption; the proof uses the geometry of nodal
domains established in Section~\ref{sec:three-dimensional-case}.

\begin{lemma}\label{lem:unconditional-sign}
Let $(M^n,g)$, $n\geq5$, be a complete, connected Riemannian manifold
without boundary, critical for $\cS^2$. Then either $S\equiv0$ or
$S<0$ everywhere.
\end{lemma}

\begin{proof}
Suppose that a connected component $\Omega$ of $\{S>0\}$ exists,
and set
\[
 u=S|_\Omega>0,\qquad
 a_n=\frac{n-4}{4(n-1)}>0,\qquad
 b_n=\frac3{4(n-1)}.
\]
The critical equations give
\[
 \Hess(u)=u\Ric-b_nu^2g,\qquad \Delta u=a_nu^2.
\]
The proof of Lemma~\ref{lem:nodal-set} applies with $b_n$ in place
of $3/8$ and with the positive side of the nodal set in place of
the negative side. Indeed, the geodesic uniqueness argument uses
only $\Hess(S)=S\Ric-b_nS^2g$; the subsequent completion and
boundary-avoidance arguments use regularity and total geodesicity
of $\{S=0\}$. Thus the intrinsic completion $X$ of $\Omega$ is a
proper geodesic space, naturally identified with $\overline\Omega$,
and minimizing segments with interior endpoints stay in $\Omega$.
The function $u$ extends continuously to $X$, with value zero on
$\partial\Omega$.

As in the three-dimensional proof, put
\[
 f=-\log u,\qquad d\mu=u\,\dv,\qquad
 L=\Delta-\langle\nabla f,\nabla\cdot\rangle.
\]
This time the identities corresponding to
\eqref{eq:three-bakry-emery} are
\begin{equation}\label{eq:positive-nodal-weighted}
 Lu=a_nu^2+\frac{|\nabla u|^2}{u}\geq a_nu^2,
 \qquad \Ric_f^1=b_nu\,g\geq0.
\end{equation}
Fix $q\in\Omega$. The distance comparison and Calabi argument
in the proof of Lemma~\ref{lem:three-weighted-package} apply in
any dimension. This time $\Ric_f^1=b_nug\geq0$, and the
corresponding radial Riccati inequality is
\[
 (L\rho)'\leq-\frac{(L\rho)^2}{n},
 \qquad \rho=d_X(q,\cdot),
\]
because the effective dimension is $n+1$. Thus the comparison
from any slightly shifted basepoint along an interior minimizing
segment gives $L\rho_\delta(x)\leq n/(\rho(x)-\delta)$.
Using the same nonincreasing function $\zeta$ as before, for
any $r>0$ the function $\eta=\zeta(\rho/r)$ has compact support
in $B_{2r}^X(q)$, equals one on $B_r^X(q)$, and at every point
$x\in\Omega$ with $\eta(x)>0$ admits a smooth lower support
$\eta_x$ satisfying
\begin{equation}\label{eq:positive-nodal-cutoff}
 \frac{|\nabla\eta_x(x)|^2}{\eta(x)}\leq\frac{C_n}{r^2},
 \qquad L\eta_x(x)\geq-\frac{C_n}{r^2}.
\end{equation}
No local bound for $u$ is needed.

Set $U=u(q)>0$ and $F=\eta u$ on $\Omega$, extended by zero
to $X\setminus\Omega$. This function is continuous and has
compact support; it reaches its maximum at some $x_*\in\Omega$,
where $F(x_*)\geq U$ and $\eta(x_*)>0$.
Let $\eta_*$ be the smooth lower support at $x_*$.
Then $F_*=u\eta_*$ has a local maximum at $x_*$, whence
$\nabla F_*(x_*)=0$ and $LF_*(x_*)\leq0$.
Using $Lu\geq a_nu^2$, we obtain at $x_*$
\begin{align*}
 0&\geq LF_*=\eta_*Lu+uL\eta_*
          +2\langle\nabla u,\nabla\eta_*\rangle\\
  &\geq u\left(a_nF(x_*)-\frac{3C_n}{r^2}\right)
   \geq u\left(a_nU-\frac{3C_n}{r^2}\right).
\end{align*}
Taking $r^2>3C_n/(a_nU)$ gives a contradiction.
Hence $S\leq0$ on $M$. Finally,
$\Delta S=a_nS^2\geq0$, so the strong maximum principle and
connectedness imply either $S\equiv0$ or $S<0$ everywhere.
\end{proof}

We may now assume $S<0$ and set $u=-S>0$. From this point onward
we impose $S\in L^p(M,g)$ for some $0<p<\infty$. Its first role
is to make $u^p\,\dv$ a finite measure. We use this measure to
bound $|\nabla u^{-1/2}|$ directly by a dimensional constant. The scalar equation
will turn the latter estimate into rigidity in the stated range
of exponents, including the upper endpoints.

\subsection{The gradient estimate}
Set $w=u^{-1/2}$ and $H=|\nabla w|$, so that
$H^2=|\nabla u|^2/(4u^3)$ is invariant under constant rescaling
of the metric. We first obtain a preliminary bound by testing a
normalized excess against the finite measure $u^p\,\dv$.
A second test, localized near the supremum of $H^2$, gives the
dimensional constant in the following proposition.

\begin{proposition}\label{prop:gradient-estimate}
Let $(M^n,g)$, $n\geq5$, be a complete critical metric of $\cS^2$,
and assume that $S<0$ and $S\in L^p(M,g)$ for some $0<p<\infty$.
Then, setting $u=-S$,
\begin{equation}\label{eq:sharp-intrinsic-gradient-estimate}
 \frac{|\nabla u|^2}{u^3}\leq4K_n,\qquad
 K_n=\max\left\{\frac1{8(n-1)},
                 \frac{(n-4)^2}{24(n-1)(n+2)}\right\}.
\end{equation}
\end{proposition}

\begin{remark}
Equivalently,
\[
 \frac{|\nabla u|^2}{u^3}\leq
 \begin{cases}
 1/[2(n-1)],&5\leq n\leq9,\\[1mm]
 (n-4)^2/[6(n-1)(n+2)],&n\geq10.
 \end{cases}
\]
\end{remark}

\begin{proof}
Set
\[
u=-S>0, \qquad w=u^{-1/2},  \qquad H=|\nabla w|.
\]
The critical system is
\begin{equation}\label{eq:negative-critical-system}
 \Ric=\frac{\Hess(u)}u-\frac3{4(n-1)}ug,\qquad
 \Delta u=-\frac{n-4}{4(n-1)}u^2.
\end{equation}
Throughout the proof set
\[
 c_n=\frac{n-4}{8(n-1)},\qquad d_n=\frac{n+2}{8(n-1)}.
\]
The critical equations and $u=w^{-2}$ give
\begin{equation}\label{eq:w-system-gradient}
 w\Delta w=3H^2+c_n,
 \qquad
 \Ric=-\frac2w\Hess(w)+\frac6{w^2}dw\otimes dw
                -\frac3{4(n-1)w^2}g.
\end{equation}
We use the symmetric diffusion operator
\[
 L_p=\Delta-\frac{2p}{w}\langle\nabla w,\nabla\cdot\rangle
     =w^{2p}\operatorname{div}(w^{-2p}\nabla\cdot),
 \qquad d\nu_p=w^{-2p}\,dV_g=u^p\,dV_g.
\]
Thus $\nu_p(M) = \int_M d\nu_p<\infty$.

Completeness provides compactly supported
Lipschitz cutoffs $\eta_r$, with
\begin{equation}\label{eq:gradient-cutoffs}
 0\leq\eta_r\leq1,\qquad \eta_r=1\text{ on }B_r(o),\qquad
 \supp\eta_r\subset B_{2r}(o),\qquad |\nabla\eta_r|\leq C/r.
\end{equation}

Bochner's formula, with \eqref{eq:w-system-gradient}, yields
\begin{equation}\label{eq:bochner-H2}
 L_p(H^2)=2|\Hess(w)|^2+\frac{6H^4-2d_nH^2}{w^2}
       +\frac{4-2p}{w}\langle\nabla w,\nabla H^2\rangle,
\end{equation}
where we used
\[
 \nabla\Delta w=\frac3w\nabla H^2-
                     \frac{3H^2+c_n}{w^2}\nabla w,
 \qquad
 2\Hess(w)(\nabla w,\nabla w)
                     =\langle\nabla w,\nabla H^2\rangle.
\]

\medskip
\textit{Step 1.}
Let
\[
 G=\frac{H^2-1/12}{w^2};
\]
the value $1/12$ is chosen in order to factorize the nonlinear term in $H^2$.
Indeed,
\[
 L_p(w^{-2})=\frac{4pH^2-2c_n}{w^4},\qquad
 \nabla(w^{-2})=-2w^{-3}\nabla w.
\]
Using these two identities, the product rule in
\eqref{eq:bochner-H2} gives, for any constant $K$,
\begin{align}
 L_p\!\left(\frac{H^2-K}{w^2}\right)
 ={}&\frac2{w^2}|\Hess(w)|^2
 -\frac{2p}{w}\left\langle\nabla w,
                \nabla\!\left(\frac{H^2-K}{w^2}\right)\right\rangle
 \nonumber\\
 &+\frac{6H^4-\frac12H^2+2c_nK}{w^4}.\label{eq:gradient-excess-general}
\end{align}
Here we used $2(c_n+d_n)=1/2$. For $K=1/12$,
$6H^4-H^2/2=6H^2(H^2-1/12)$, so
\begin{equation}\label{eq:gradient-excess-exact}
 L_pG=\frac2{w^2}|\Hess(w)|^2
       -\frac{2p}{w}\langle\nabla w,\nabla G\rangle
       +6\frac{H^2}{w^2}G+\frac{c_n}{6w^4}.
\end{equation}
Let $G_+=\max\{G,0\}$. Choose an integer $k\geq4$ with
$k-2\geq p^2$ and put $\eta=\eta_r$.
We multiply \eqref{eq:gradient-excess-exact} directly by
$\eta^{2k}G_+^{k-2}$ and integrate with respect to $d\nu_p$.
This is an admissible compactly supported Lipschitz test: since
$k\geq4$, the function $t\mapsto(t_+)^{k-2}$ is continuously
differentiable and
\[
 \nabla(G_+^{k-2})=(k-2)G_+^{k-3}\nabla G.
\]
The test vanishes where $G\leq0$, while $G=G_+$ where it is positive;
note that the terms involving $|\Hess(w)|^2$ and $c_n$ in
\eqref{eq:gradient-excess-exact} are nonnegative and may be discarded.
Let
\[
 A=\int_M\eta^{2k}\frac{H^2}{w^2}G_+^{k-1}\,d\nu_p,\qquad
 B=\int_M\eta^{2k}G_+^{k-3}|\nabla G_+|^2\,d\nu_p;
\]
integration by parts, using $|\nabla w|=H$ and the Cauchy--Schwarz inequality, gives
\begin{align*}
 6A+(k-2)B\leq{}&
 2k\int_M\eta^{2k-1}G_+^{k-2}|\nabla\eta|\,|\nabla G_+|\,d\nu_p\\
 &+2p\int_M\eta^{2k}\frac{H}{w}G_+^{k-2}|\nabla G_+|\,d\nu_p.
\end{align*}
Young's inequality bounds the first term by
\[
 \frac{k-2}{2}B+
 \frac{2k^2}{k-2}\int_M\eta^{2k-2}G_+^{k-1}|\nabla\eta|^2\,d\nu_p
\]
and the second by
\[
 \frac{k-2}{2}B+\frac{2p^2}{k-2}A
 \leq\frac{k-2}{2}B+2A.
\]
Consequently,
\begin{equation}\label{eq:gradient-excess-test}
 4A\leq\frac{2k^2}{k-2}
               \int_M\eta^{2k-2}G_+^{k-1}|\nabla\eta|^2\,d\nu_p.
\end{equation}
Set $I_r=\int_M\eta_r^{2k}G_+^k\,d\nu_p$.
Since $\frac{H^2}{w^2}\geq G_+$, we have
$I_r\leq A$. H\"older's inequality with conjugate exponents
$k/(k-1)$ and $k$ yields
\[
 I_r\leq C_k I_r^{(k-1)/k}
          \left(\int_M|\nabla\eta_r|^{2k}\,d\nu_p\right)^{1/k}.
\]
If $I_r>0$ we divide and raise to the power $k$; if $I_r=0$ the
resulting inequality is automatic. Thus
\begin{equation}\label{eq:gradient-excess-limit}
 I_r\leq C_{k,p}\,r^{-2k}\nu_p(M)\longrightarrow0.
\end{equation}
On any fixed ball $\eta_r=1$ for all sufficiently large $r$;
it follows that $G_+=0$ there, and hence everywhere. We have proved
\begin{equation}\label{eq:gradient-preliminary}
 H^2\leq\frac1{12}\quad\text{on }M.
\end{equation}

\medskip
\textit{Step 2.}
At a point with $H>0$, choose an orthonormal frame with
$e_1=\nabla w/H$. The sum of the remaining diagonal Hessian
entries is $\Delta w-\Hess(w)(e_1,e_1)$. Cauchy--Schwarz on these
$n-1$ entries, together with nonnegativity of the omitted squares, gives
\begin{align*}
 |\Hess(w)|^2
 &\geq \Hess(w)(e_1,e_1)^2
       +\frac{\bigl(\Delta w-\Hess(w)(e_1,e_1)\bigr)^2}{n-1}\\
 &=\frac{(\Delta w)^2-2\Delta w\,\Hess(w)(e_1,e_1)}{n-1}
       +\frac{n}{n-1}\Hess(w)(e_1,e_1)^2\\
 &\geq\frac{(\Delta w)^2-2\Delta w\,\Hess(w)(e_1,e_1)}{n-1}.
\end{align*}
Since
\[
 \Delta w=\frac{3H^2+c_n}{w},\qquad
 2H^2\Hess(w)(e_1,e_1)=\langle\nabla w,\nabla H^2\rangle,
\]
\eqref{eq:bochner-H2} implies
\begin{equation}\label{eq:sharp-scalar-bochner}
 L_p(H^2)\geq\frac{F_n(H^2)}{w^2}
   +\frac{b_p(H^2)}w\langle\nabla w,\nabla H^2\rangle,
\end{equation}
where
\begin{align}
 F_n(s)&=\frac{2(3s+c_n)^2}{n-1}+6s^2-2d_ns
       =\frac{6(n+2)}{n-1}(s-\kappa_1)(s-\kappa_2),\label{eq:sharp-quadratic-factorization}\\
 b_p(s)&=4-2p-\frac{2(3+c_n/s)}{n-1},\qquad
 \kappa_1=\frac1{8(n-1)},\quad
 \kappa_2=\frac{(n-4)^2}{24(n-1)(n+2)}.\nonumber
\end{align}
In particular, $F_n(s)>0$ for $s>K_n$.

By \eqref{eq:gradient-preliminary}, the number $M_0=\sup_M H^2$ is finite.
Suppose that $M_0>K_n$ and put $m_0=(M_0+K_n)/2$.
On the compact interval $[m_0,M_0]$, choose constants
\[
 a_0=\min_{[m_0,M_0]}F_n>0,
 \qquad b_0=\max_{[m_0,M_0]}|b_p(s)|\sqrt{s}<\infty.
\]
Take $0<\varepsilon<M_0-m_0$ so small that
$b_0^2\varepsilon\leq a_0/2$, and define
\[
 K=M_0-\varepsilon,\qquad Z=(H^2-K)_+.
\]
Then $0\leq Z\leq\varepsilon$ and $Z$ is positive somewhere, by
the definition of the supremum. On $\{Z>0\}$, inequality
\eqref{eq:sharp-scalar-bochner} gives
\[
 L_p(H^2)\geq\frac{a_0}{w^2}-\frac{b_0}{w}|\nabla(H^2)|.
\]
We test this inequality directly with $\eta_r^2Z$.

The support of this Lipschitz test lies in the compact set
$\supp\eta_r\cap\{H^2\geq K\}\subset\{H>0\}$, because $K>0$.
Thus the coefficient $b_p(H^2)$ is bounded near its support.

The first-order chain rule gives, almost everywhere,
\[
 \langle\nabla(H^2),\nabla Z\rangle=|\nabla Z|^2,
 \qquad Z\nabla(H^2)=Z\nabla Z.
\]
Integration by parts applied to the smooth function $H^2$ therefore
gives, writing $\eta=\eta_r$,
\begin{align*}
 \int_M\eta^2|\nabla Z|^2\,d\nu_p
 +a_0\int_M\eta^2\frac{Z}{w^2}\,d\nu_p
 \leq{}&2\int_M\eta Z|\nabla\eta|\,|\nabla Z|\,d\nu_p\\
 &+b_0\int_M\eta^2\frac{Z}{w}|\nabla Z|\,d\nu_p.
\end{align*}
Using again Young's inequality, the first term on the right is at most
$\frac14\int\eta^2|\nabla Z|^2\,d\nu_p+4\int Z^2|\nabla\eta|^2\,d\nu_p$, while
the second is at most
$\frac14\int\eta^2|\nabla Z|^2\,d\nu_p+b_0^2\int\eta^2Z^2/w^2\,d\nu_p$.
Using $Z^2\leq\varepsilon Z$ and the choice of $\varepsilon$, we obtain
\begin{equation}\label{eq:gradient-near-supremum}
 \frac12\int_M\eta_r^2|\nabla Z|^2\,d\nu_p
 +\frac{a_0}{2}\int_M\eta_r^2\frac{Z}{w^2}\,d\nu_p
 \leq4\int_M Z^2|\nabla\eta_r|^2\,d\nu_p
 \leq\frac{C\varepsilon^2}{r^2}\nu_p(M).
\end{equation}
Let $r\to\infty$ while keeping a fixed ball on the left.
We conclude that $Z/w^2$ has zero integral on that ball.
Since $w>0$ and $Z$ is continuous, $Z=0$ everywhere, contradicting
its positivity somewhere. Thus $M_0\leq K_n$, proving
\eqref{eq:sharp-intrinsic-gradient-estimate}.
\end{proof}

\begin{remark}
\label{rem:gradient-sharpness-forward}
The constant $K_n$ in Proposition~\ref{prop:gradient-estimate}
is optimal in dimension $n=6$ and in every dimension $n\geq10$.
Indeed, the complete examples constructed in
Section~\ref{sec:direct-examples} have scalar curvature in $L^p$
for suitable $p>1$ and satisfy
\[
 \sup_M|\nabla(-S)^{-1/2}|^2=K_n;
\]
see Remark~\ref{rem:sharp-gradient-examples} for the explicit
formulas. In these examples the supremum is approached at infinity,
not attained at a finite point.
\end{remark}

\subsection{Rigidity}
The gradient estimate is valid for every integrability exponent
$p>0$. The restriction $p\leq P_n$ enters in the following argument,
where the scalar equation is integrated against powers of $u$.

\begin{proof}[Proof of Theorem \ref{thm:rigidity-intro}]
The compact case follows immediately by integrating
\[
\Delta S=\frac{n-4}{4(n-1)}S^2.
\]

Assume therefore that $M$ is noncompact. By
Lemma~\ref{lem:unconditional-sign}, either $S\equiv0$ or $S<0$
everywhere. In the latter case set $u=-S>0$ and argue by
contradiction.

Proposition~\ref{prop:gradient-estimate} gives
\begin{equation}\label{eq:gradient-bound-final}
\frac{|\nabla u|^2}{u^3}
\leq
4K_n,
\qquad
K_n
=
\max\left\{
\frac{1}{8(n-1)},
\frac{(n-4)^2}{24(n-1)(n+2)}
\right\}.
\end{equation}
Since $w=u^{-1/2}$,
\[
|\nabla w|\leq\sqrt{K_n}.
\]
Thus, fixing $o\in M$ and integrating along minimizing geodesics,
\begin{equation}\label{eq:w-growth-final}
u(x)^{-1/2}
\leq
u(o)^{-1/2}
+
\sqrt{K_n}\,d(x,o).
\end{equation}

Let $\eta_s$ be a standard radial cutoff with $0\leq\eta_s\leq1$,
equal to one on $B_s(o)$, supported in $B_{2s}(o)$, and with
$|\nabla\eta_s|\leq C/s$, as in \eqref{eq:gradient-cutoffs}.

Multiplying

\[
-\Delta u
=
\frac{n-4}{4(n-1)}u^2
\]
by $u^{p-2}\eta_s^2$ and integrating by parts gives
\begin{align}
\frac{n-4}{4(n-1)}
\int_Mu^p\eta_s^2\,\dv
={}&
(p-2)
\int_M
\frac{|\nabla u|^2}{u^3}
u^p\eta_s^2\,\dv
\nonumber\\
&+
2\int_M
u^{p-2}\eta_s
\langle\nabla u,\nabla\eta_s\rangle\,\dv.
\label{eq:final-integral-identity}
\end{align}
The last term in \eqref{eq:final-integral-identity} tends to zero. Indeed,
\eqref{eq:gradient-bound-final} gives
\[
|\nabla u|
\leq
2\sqrt{K_n}\,u^{3/2},
\]
while \eqref{eq:w-growth-final} gives
\[
u^{-1/2}\leq Cs
\]
on $B_{2s}(o)\setminus B_s(o)$ for all sufficiently large $s$.
Consequently,
\begin{align*}
\left|
2\int_M
u^{p-2}\eta_s
\langle\nabla u,\nabla\eta_s\rangle\,\dv
\right|
&\leq
\frac{C}{s}
\int_{B_{2s}(o)\setminus B_s(o)}
u^{p-\frac12}\,\dv\\
&\leq
C
\int_{B_{2s}(o)\setminus B_s(o)}
u^p\,\dv
\longrightarrow0.
\end{align*}

Here $u^{p-1/2}=w u^p$, so the estimate uses exactly the
assumed integrability of $u^p$, together with the linear growth
bound for $w$. No volume growth condition or inclusion between
different global $L^p$ spaces is required.

If $0<p\leq2$, the first term on the right-hand side of
\eqref{eq:final-integral-identity} is nonpositive. Letting
$s\to\infty$ therefore gives
\[
\frac{n-4}{4(n-1)}
\int_Mu^p\,\dv\leq0,
\]
which is impossible. The global integrability needed to pass to the limit follows
from the gradient estimate:
\[
 0\leq u^{p-3}|\nabla u|^2\leq4K_nu^p\in L^1(M).
\]
This holds for the original exponent $p$, including $0<p<2$.

Suppose now that $p>2$. By \eqref{eq:gradient-bound-final},
\[
\left[
\frac{n-4}{4(n-1)}
-
4(p-2)K_n
\right]
\int_Mu^p\eta_s^2\,\dv
\leq
C
\int_{B_{2s}(o)\setminus B_s(o)}
u^p\,\dv.
\]

The threshold is determined by the balance between the two coefficients:
\[
 \frac{n-4}{4(n-1)}-4(p-2)K_n=4K_n(P_n-p),\qquad
 P_n=2+\frac{n-4}{16(n-1)K_n}.
\]
Substitution of the two branches of $K_n$ yields
\[
 P_n=\begin{cases}
 2+(n-4)/2=n/2,&5\leq n\leq9,\\[1mm]
 2+3(n+2)/[2(n-4)]=(7n-10)/[2(n-4)],&n\geq10.
 \end{cases}
\]
For $p<P_n$, the coefficient is positive, and letting $s\to\infty$
contradicts $u>0$ and $\int_Mu^p\,\dv<\infty$.

It remains to consider the endpoint. There
\[
\frac{n-4}{4(n-1)}
=
4(p-2)K_n.
\]
Letting $s\to\infty$ in the exact identity
\eqref{eq:final-integral-identity}, using the vanishing of the

boundary term and dominated convergence, gives the following identity.
The domination follows from
$0\leq|\nabla u|^2/u^3\leq4K_n$ and $u^p\in L^1$:

\[
\int_M
\left(
4K_n-\frac{|\nabla u|^2}{u^3}
\right)
u^p\,\dv
=0.
\]
The integrand is continuous and nonnegative, hence
\[
\frac{|\nabla u|^2}{u^3}
\equiv4K_n.
\]
Equivalently,
\[
|\nabla w|^2\equiv K_n>0.
\]

The vector field
\[
 Z=-\frac{\nabla w}{\sqrt{K_n}}
\]
is smooth and has unit length. Its integral curves exist for all
positive time. Indeed, on any finite time interval they stay in a
closed bounded ball, which is compact by completeness; the smooth
ordinary differential equation can therefore be extended beyond
any finite endpoint of its interval of existence. Along such a
curve $\gamma$,

\[
\frac{\mathrm{d}}{\mathrm{d}t}w(\gamma(t))
=
-|\nabla w|
=
-\sqrt{K_n}.
\]
Therefore
\[
w(\gamma(t))
=
w(\gamma(0))-\sqrt{K_n}\,t,
\]
so $w$ reaches zero in finite positive time. This contradicts
$w=u^{-1/2}>0$.

Hence the alternative $S<0$ is impossible, and therefore
$S\equiv0$.
\end{proof}

\begin{remark}

Note that, since
\[
 \frac{(n-4)^2}{24(n-1)(n+2)}-\frac1{8(n-1)}
                         =\frac{n-10}{24(n+2)},
\]
at $n=10$ both gradient constants agree and the endpoint is
$P_{10}=5$.
The examples in Theorem~\ref{thm:direct-examples} below prove that
$P_n$ is optimal for every $n\geq10$: their scalar curvature belongs
to $L^p$ exactly for $p>P_n$.

\end{remark}

  The sign lemma also yields the following finite-volume rigidity statement.
\begin{corollary}
Let $(M^n, g)$ be a complete, connected, critical metric for $\cS^2$, with $n\geq 5$. If $\operatorname{Vol}_g(M)<\infty$, then $S_g\equiv 0$.
\end{corollary}
\begin{proof}
By Lemma~\ref{lem:unconditional-sign}, it suffices to exclude
$S<0$. Set $u=-S>0$ and $a_n=\frac{n-4}{4(n-1)}$. Testing
$-\Delta u=a_nu^2$ with $\eta^2/u$, where $\eta$ is a
compactly supported cutoff, gives
\[
 a_n\int_M\eta^2u\,\dv
 =-\int_M\eta^2|\nabla\log u|^2\,\dv
 +2\int_M\eta\langle\nabla\log u,\nabla\eta\rangle\,\dv
 \leq\int_M|\nabla\eta|^2\,\dv.
\]
Choosing standard cutoffs $\eta_R$ with $\eta_R=1$ on $B_R$, $\eta_R=0$ on $M\setminus B_{2R}$
and $|\nabla\eta_R|\leq C/R$, the right-hand side tends to
zero because $\vol(M)<\infty$.

For every  ball $B_{R_0}$ and $R\geq R_0$ we obtain
\[
a_n\int_{B_{R_0}}u\,\dv\leq C R^{-2}\vol(M).
\]
Letting $R\to\infty$ gives $u\equiv 0$ on that ball and then on $M$,
a contradiction.

\end{proof}

\section{Explicit complete examples in higher dimensions}
\label{sec:direct-examples}

In this section we construct an explicit family showing that completeness,
bounded scalar curvature, and its pointwise decay to zero at infinity do
not imply scalar-flatness. The metrics are defined on the $n$-dimensional product $\mathbb{R}^2 \times \mathbb{T}^{n-2}$, where $\mathbb{T}^{n-2}$ is a flat torus of dimension $n-2$ and $n\geq 5$.

For the sake of completeness, we first recall a few basic facts about moving frames (see for instance \cite{AMR} or \cite{CMBook}).

Fix a local orthonormal coframe $\set{\theta^i}$ for $g$, with
$i=1,\ldots,n$, and let $\set{e_i}$ be its dual frame. We follow the Einstein summation convention: repeated
indices in formulas are summed from $1$ to $n$,
unless otherwise specified. The corresponding
\emph{Levi-Civita connection forms}
$\set{\theta^i_j}$  are the $1$-forms uniquely defined by the
requirements
\begin{align}
  &d\theta^i = -\theta^i_j \wedge
  \theta^j \quad \text{(first structure equations)}, \label{1_firstStructureEq} \\
  &\theta^i_j + \theta^j_i = 0. \label{1_skewsymmConnForm}
\end{align}
The \emph{curvature forms} $\set{\Theta^i_j}$ associated to the
connection are the $2$-forms defined via the \emph{second structure
equations}
\begin{equation}\label{1_secondStructureEq}
  {d\theta^i_j = -\theta^i_k \wedge \theta^k_j + \Theta^i_j.}
\end{equation}
They are skew-symmetric (i.e. $\Theta^i_j + \Theta^j_i = 0$) and they can be written as
\begin{equation}\label{1_def_forme_curvatura}
  \Theta^i_j = \frac{1}{2}R^i_{jkt}\theta^k \wedge \theta^t = \sum_{k<t}R^i_{jkt}\theta^k \wedge \theta^t,
\end{equation}
where $R^i_{jkt}$ are the components of the Riemann curvature
tensor in its $(1,3)$ version. Since we are working on an orthonormal coframe, the coefficients of the $(0, 4)$ version of the Riemann tensor are simply given by
\[
R_{ijkt} = R^i_{jkt},
\]
while the coefficients of the Ricci tensor are obtained by tracing:
\[
R_{ik} = \sum_{t=1}^n R_{itkt}.
\]
Similarly, the scalar curvature $S$ is the trace of the Ricci tensor, that is
\[
S = \sum_{t=1}^n R_{tt}.
\]

 We use commas for covariant derivatives of tensor components,
e.g. $R_{ij,k}$. For a scalar function we instead write subscripts
without commas: in particular, $S_i$ and $S_{ij}$ denote the
components of $dS$ and of the Hessian of $S$, respectively.
The covariant derivative of a vector field $X \in \mathfrak{X}(M)$ is defined
 by
\[
dX^i+X^j\theta^i_j=X^i_{,k}\theta^k,
\]
while the {covariant derivative of a $1$-form} $\omega$ is
defined by
\[
d\omega_i-\omega_j\theta^j_i=\omega_{i,k}\theta^k.
\]

The covariant derivative of $X$ defines an endomorphism of $TM$
whose matrix in the orthonormal frame is $(X^i_{,j})$.
Its trace is the divergence:
\begin{equation}\label{1_divergence}
 \operatorname{div}X=X^i_{,i}.
\end{equation}

For a smooth function $f$ we can write
\begin{equation}\label{DifferentialComponents}
df = f_i\theta^i,
\end{equation}
for smooth coefficients $f_i$ defined on the coframe domain. The Hessian of $f$, $\Hess(f)$, is the $(0, 2)$-tensor defined as
\begin{equation}
  \Hess(f) = f_{ij}\theta^j\otimes\theta^i,
\end{equation}
with
\begin{equation}\label{HessianComponents}
f_{ij}\theta^j = df_i - f_t\theta_i^t.
\end{equation}
Note that
\[
f_{ij}=f_{ji}.
\]

The Laplacian of $f$, $\Delta f$,  is the trace of the
Hessian, in other words
\[
\Delta f = \operatorname{tr}(\Hess(f)) = f_{ii}.
\]

We are now ready to state the following theorem, which immediately implies Theorem~\ref{thm:examples-intro} in the Introduction.

\begin{theorem}
\label{thm:direct-examples}
Let $n\geq5$, let $g_0$ be a flat metric on $\T^{n-2}$, and set
\begin{equation}\label{eq:direct-parameters}
 a=-\frac{n+14}{6(n+2)},\qquad
 b=-\frac{n-2}{2(n+2)},\qquad
 c=\frac2{n+2},\qquad
 \lambda^2=\frac{32(n-1)}{3(n+2)},\quad \lambda>0.
\end{equation}
In standard polar coordinates $(r,\theta)$ on $\R^2$, with
$\theta\in\R/(2\pi\mathbb Z)$, define
\begin{equation}\label{eq:direct-metric}
 \begin{split}
 g_n&=A(r)^2\,dr^2+B(r)^2\,d\theta^2+C(r)^2g_0,\\
 A(r)&=\lambda(1+r^2)^a,\qquad
 B(r)=\lambda r(1+r^2)^b,\qquad
 C(r)=(1+r^2)^c.
 \end{split}
\end{equation}
Then $g_n$ is a smooth complete critical metric of $\cS^2$ on
$\R^2\times\T^{n-2}$, with
\begin{equation}\label{eq:direct-scalar}
 S_{g_n}=-(1+r^2)^{-(2a+1)}
 =-(1+r^2)^{-\frac{2(n-4)}{3(n+2)}}.
\end{equation}
In particular, $g_n$ is neither scalar-flat nor Einstein, its scalar curvature
is bounded, and $S_{g_n}\to0$ at infinity. Moreover for every $p>0$
\begin{equation}\label{eq:direct-integrability}
 S_{g_n}\in L^p(M,g_n)\quad\Longleftrightarrow\quad p>q_n,
\end{equation}
where
\[
 q_n=\frac{7n-10}{2(n-4)}.
\]
Writing $V_0=\vol_{g_0}(\T^{n-2})$, one has, for $p>q_n$,
\begin{equation}\label{eq:direct-integral}
 \int_M|S_{g_n}|^p\,\dv_{g_n}
 =\frac{16\pi(n-1)V_0}{(n-4)(p-q_n)}.
\end{equation}
\end{theorem}

\begin{proof}
Set $\Sigma=\{0\}\times\T^{n-2}$. Throughout the proof, Greek indices
$\alpha,\beta,\gamma\in\{3,\ldots,n\}$ denote directions tangent to the torus.

\textbf{Step 1.} First we need to choose a suitable frame, and then to compute all the relevant geometric quantities.
For $r>0$, we choose a local, parallel $g_0$-orthonormal coframe
$\{\eta^\alpha\}_{\alpha=3}^n$ on the torus (note that such a choice is possible
because $g_0$ is flat, and $d\eta^\alpha=0$); thus, an orthonormal coframe
for $g_n$ is given by
\[
 \theta^1=A\,dr,\qquad \theta^2=B\,d\theta,\qquad
 \theta^\alpha=C\eta^\alpha\quad(\alpha=3,\ldots,n).
\]
For the moment, $A,B,C$ are arbitrary positive smooth functions
on $r>0$; primes denote derivatives with respect to $r$. The Levi-Civita connection forms are given by
\begin{align*}
  \theta^2_1 &= -\theta^1_2 = \frac{B'}{AB}\theta^2, \\
  \theta^\alpha_1 &= -\theta^1_\alpha = \frac{C'}{AC}\theta^\alpha,\\
  \theta^\alpha_2 &= -\theta^2_\alpha = 0,\\
  \theta^\alpha_\beta &= -\theta^\beta_\alpha = 0 \,\, (\alpha, \beta = 3, \ldots, n).
\end{align*}

The independent nonzero components of the Riemann tensor are
\begin{align*}
 R_{1212}&=\frac1{A^2B}\left(\frac{A'}A B'-B''\right),\\
 R_{1\alpha1\alpha}&=\frac1{A^2C}\left(\frac{A'}A C'-C''\right),\\
 R_{2\alpha2\alpha}&=-\frac1{A^2}\frac{B'C'}{BC},\\
 R_{\alpha\beta\alpha\beta}&=-\frac{(C')^2}{A^2C^2}\qquad(\alpha\neq\beta).
\end{align*}
All other components vanish, apart from those determined by the
curvature symmetries. There is no summation over $\alpha$ or $\beta$
in these component formulas or in the diagonal formulas below.
Tracing, we find that the Ricci tensor is diagonal, with components

\begin{align*}
 R_{11}
 &=-\frac1{A^2}\left[
 \frac{B''}{B}-\frac{A'B'}{AB}
 +(n-2)\left(\frac{C''}{C}-\frac{A'C'}{AC}\right)\right],\\
 R_{22}
 &=-\frac1{A^2}\left[
 \frac{B''}{B}-\frac{A'B'}{AB}
 +(n-2)\frac{B'C'}{BC}\right],\\
 R_{\alpha\alpha}
 &=-\frac1{A^2}\left[
 \frac{C''}{C}-\frac{A'C'}{AC}+\frac{B'C'}{BC}
 +(n-3)\left(\frac{C'}C\right)^2\right].
\end{align*}
Finally, tracing again we obtain the expression for the scalar curvature:

\begin{align*}
  S &= \frac{2}{A^2B}\pa{\frac{A'}{A}B'-B''}+ \frac{2(n-2)}{A^2C}\pa{\frac{A'}{A}C'-C''} - \frac{2(n-2)}{A^2} \frac{B'C'}{BC} - \frac{(n-2)(n-3)}{(AC)^2}\pa{C'}^2.
\end{align*}
Since $S$ depends only on $r$, equation \eqref{HessianComponents} shows that its Hessian is also diagonal, with
\begin{align*}
  S_{11} &= \frac{1}{A^2}\pa{S''-\frac{A'}{A}S'},\\
S_{22} &= \frac{1}{A^2}\frac{B'}{B}S',\\
S_{\alpha\alpha} &= \frac{1}{A^2}\frac{C'}{C}S'.
\end{align*}

\textbf{Step 2.} We now specialize the choice of $A, B$ and $C$, proving the smoothness of the corresponding metric and then imposing the critical equation.
We look for a critical metric within the family
\[
 A=\lambda(1+r^2)^a,\qquad
 B=\lambda r(1+r^2)^b,\qquad
 C=(1+r^2)^c,
\]
where $a,b,c\in\R$ and $\lambda>0$ are initially
unknown.

\begingroup
\noindent We prove that the metric extends smoothly
across $\Sigma$ for every
$a,b,c\in\mathbb R$ and $\lambda>0$. Indeed, factoring out
$A(r)^2$, we can write
\[
 \begin{split}
 g_n=\lambda^2(1+r^2)^{2a}\Bigl[dr^2
       +r^2(1+r^2)^{2(b-a)}d\theta^2
       +\lambda^{-2}(1+r^2)^{2(c-a)}g_0\Bigr].
 \end{split}
\]
The two-dimensional part of the metric in brackets is
$dr^2+\varphi(r)^2d\theta^2$, where
\[
 \varphi(r)=r(1+r^2)^{b-a}.
\]
This function is smooth and odd, with
$\varphi(0)=0$ and $\varphi'(0)=1$. Since $\theta$ has
period $2\pi$, the standard endpoint smoothness criterion
implies that this two-dimensional metric extends smoothly
across the origin; see
\cite[Chapter~1, Sections~3--4]{petersen}.
Adding the torus term and multiplying by the conformal factor
therefore gives a smooth Riemannian metric on the whole
$\mathbb R^2\times\mathbb T^{n-2}$.

\endgroup

We now impose the critical equation and the normalization
$S|_\Sigma=-1$; these conditions will determine all four parameters.
Put $x=r^2$. The logarithmic derivatives are
\[
 \frac{A'}A=\frac{2a r}{1+x},\qquad
 \frac{B'}B=\frac1r+\frac{2b r}{1+x},\qquad
 \frac{C'}C=\frac{2c r}{1+x}.
\]
To simplify the writing, define the following coefficients, all independent of $r$:
\[
\begin{aligned}
 U&=3b-a+(n-2)c,\\
 V_1&=U+2b(b-a-1)
             +2(n-2)c(c-a-1),\\
 V_2&=U+2b(b-a-1)+2(n-2)bc,\\
  J&=b-a+(n-2)c.
\end{aligned}
\]
The Ricci components therefore become
\begin{equation}\label{eq:direct-Ricci}
\begin{aligned}
 R_{11}&=-\frac{2(U+V_1x)}{\lambda^2(1+x)^{2a+2}},\qquad
 R_{22}=-\frac{2(U+V_2x)}{\lambda^2(1+x)^{2a+2}},\\
 R_{\alpha\alpha}&=-\frac{4c(1+Jx)}{\lambda^2(1+x)^{2a+2}}.
\end{aligned}
\end{equation}
Taking their trace gives
\begin{equation}\label{eq:direct-scalar-polynomial}
 S=\frac4{\lambda^2}(1+x)^{-2a-2}(P_0+P_2x),
\end{equation}
where
\[
\begin{aligned}
 P_0&=a-3b-2(n-2)c,\\
 P_2&=(a-b)\bigl(1+2b+2(n-2)c\bigr)
                       -(n-2)(n-1)c^2.
\end{aligned}
\]
In particular, the prescribed value at $\Sigma$ requires
\begin{equation}\label{eq:direct-normalization-forced}
 P_0<0,\qquad \lambda^2=-4P_0.
\end{equation}
The first derivative of the scalar curvature can be written as
\[
 S'=\frac{8r}{\lambda^2}(1+x)^{-2a-3}
       \bigl(P_1-(2a+1)P_2x\bigr),\qquad
 \text{where }\,\, P_1=P_2-2(a+1)P_0.
\]
We first use the ``torus component'' of the critical equation, that is
\[
 SR_{\alpha\alpha}-S_{\alpha\alpha}=\frac3{4(n-1)}S^2.
\]
The Hessian formula above gives
\[
 S_{\alpha\alpha}
 =\frac{16c x}{\lambda^4(1+x)^{4a+4}}
       \bigl(P_1-(2a+1)P_2x\bigr).
\]
If $c=0$, both $R_{\alpha\alpha}$ and $S_{\alpha\alpha}$ vanish,
so the critical equation would force $S=0$, contrary to the
normalization; hence, necessarily, $c\neq0$. Cancelling the common nonzero
factor, the torus equation gives the polynomial identity
\begin{equation}\label{eq:direct-torus-polynomial}
\begin{split}
0={}&(P_0+P_2x)(1+Jx)
       +x\bigl(P_1-(2a+1)P_2x\bigr)\\
 &+\frac3{4c(n-1)}(P_0+P_2x)^2.
\end{split}
\end{equation}
Its constant coefficient, since $P_0\neq0$, gives
\[
 P_0=-\frac43(n-1)c;
\]
in particular, $c>0$. Using this relation, the coefficient of
$x$ becomes $P_0J+P_1-P_2=0$, hence $J=2a+2$;
the coefficient of $x^2$ then becomes
$P_2(1-P_2/P_0)=0$. We have obtained the three necessary conditions
\begin{equation}\label{eq:direct-torus-conditions}
 P_0=-\frac43(n-1)c,\qquad
 b-3a-2+(n-2)c=0,\qquad
 P_2=0\ \text{or}\ P_2=P_0.
\end{equation}
Conversely, these conditions make every coefficient in
\eqref{eq:direct-torus-polynomial} vanish.

The first two conditions are linear in $a,b$ and give
\begin{equation}\label{eq:direct-linear-parameters}
 a=\frac{7n-10}{24}c-\frac34,\qquad
 b=\frac{6-n}{8}c-\frac14.
\end{equation}
Substitution in the expression for $P_2$ yields
\[
 48P_2=-(13n-22)(n+2)c^2-32(n-1)c-12<0.
\]
Here $n\geq5$ and $c>0$. The branch $P_2=0$ is therefore
impossible, and necessarily
\[
 P_2=P_0,\qquad P_1=-(2a+1)P_0.
\]

The remaining relation $P_2=P_0$ is equivalent to
\begin{equation}\label{eq:direct-two-branches}
 \bigl((n+2)c-2\bigr)\bigl((13n-22)c-6\bigr)=0.
\end{equation}
Thus the torus component alone leaves two possible values:
\[
 c=\frac2{n+2}\qquad\text{or}\qquad c=\frac6{13n-22}.
\]
The radial component will distinguish between these two possibilities.

In particular, the scalar curvature has simplified to a single power:
\[
 S=\frac{4P_0}{\lambda^2}(1+x)^{-k},\qquad k=2a+1.
\]

\begingroup
We next evaluate the radial component of the critical equation
along $\Sigma$. Although the polar orthonormal frame is not defined
at $r=0$, the smoothness   established above shows
that $\Ric$ and $\Hess(S)$ are smooth there.
By rotational invariance, their radial components have limits independent of
the direction of approach. Writing $R_{11}(0)$ and $S_{11}(0)$ for
these common normal-directional limits, the formulas above give
\endgroup
\[
 R_{11}(0)=-\frac{2U}{\lambda^2},\qquad
 S_{11}(0)=\frac{8P_1}{\lambda^4}.
\]
Thus $SR_{11}-S_{11}=3S^2/[4(n-1)]$ at $r=0$ is equivalent to
\[
 P_0U+P_1+\frac{3P_0^2}{2(n-1)}=0,
\]
or, after the torus relations,
\[
 0=U-k-2c
   =3b-3a-1+(n-4)c
   =\frac12-\frac{n+2}{4}c.
\]
Thus $c=2/(n+2)$. Equations
\eqref{eq:direct-linear-parameters} and
\eqref{eq:direct-normalization-forced} now give exactly
\eqref{eq:direct-parameters}. Thus these constants are forced within
the chosen family by criticality and $S|_\Sigma=-1$.

\begingroup
\medskip
\textbf{Step 3.} We now verify the validity of the full critical equation.
So far, the parameters have been determined using the normalization
$S|_\Sigma=-1$, the torus component of the critical equation for
$r>0$, and its radial component at $r=0$. These relations are
necessary, but they do not
by themselves ensure that the tensor equation is satisfied for every
$r$. To establish that the resulting metric is indeed critical,
we now check all its independent components on $r>0$.
\endgroup
Substitution of the parameters just found gives
\[
 P_0=P_2=-\frac{8(n-1)}{3(n+2)},\qquad
 J=2a+2,\qquad P_1=-(2a+1)P_0.
\]
Thus all three relations in \eqref{eq:direct-torus-conditions}
hold, so the torus equation is satisfied for every $r>0$.
To check the two remaining components, direct differentiation gives
\[
\begin{aligned}
 S_{11}&=-\frac{8kP_0}{\lambda^4(1+x)^{4a+3}}(1-3kx),\\
 S_{22}&=-\frac{8kP_0}{\lambda^4(1+x)^{4a+3}}
                                  \bigl(1+(1+2b)x\bigr).
\end{aligned}
\]
The constant coefficients of both critical equations are
$U-k-2c$, which is already zero. Their coefficients of $x$
are, respectively, $V_1+3k^2-2c$ and
$V_2-k(1+2b)-2c$. To evaluate these, write $D=n+2$.
The parameters just found give
\[
\begin{gathered}
 c=\frac2D,\quad
 k=\frac{2(D-6)}{3D},\quad U=\frac23,\quad
 b-a-1=-\frac43+\frac4D,\quad
 c-a-1=-\frac56+\frac4D,\\
 V_1=\frac4D-\frac{4(D-6)^2}{3D^2}=2c-3k^2,\qquad
 V_2=\frac4D+\frac{8(D-6)}{3D^2}
                         =2c+k(1+2b).
\end{gathered}
\]
All coefficients therefore vanish. The torus component was already
settled by \eqref{eq:direct-torus-conditions}, and all mixed
components vanish. We have proved
\begin{equation}\label{eq:direct-tensor}
 S_{g_n}\Ric_{g_n}-\Hess_{g_n}(S_{g_n})
                 =\frac3{4(n-1)}S_{g_n}^2g_n
\end{equation}
on $r>0$. Moreover, $4P_0/\lambda^2=-1$ gives
\eqref{eq:direct-scalar}.
\begingroup
Since $g_n$ is smooth across $\Sigma$, as proved above, both sides
of \eqref{eq:direct-tensor} are smooth on the whole manifold;
thus the identity extends from $r>0$ to $\Sigma$ by continuity.
\endgroup
Tracing \eqref{eq:direct-tensor} yields
the scalar equation \eqref{eq:EL-trace}, so the full critical system
holds.

\textbf{Step 4.} We now check the completeness of $g_n$: we prove that every closed metric ball is compact, and then apply
the Hopf--Rinow theorem. The key point is that reaching large values of the
Euclidean radius $r$ requires arbitrarily large $g_n$-length,
even for curves that also move in the angular or torus directions.

Define
\begin{equation}\label{eq:direct-radial-length}
 \rho(r)=\int_0^r A(t)\,dt
        =\lambda\int_0^r(1+t^2)^a\,dt.
\end{equation}
This is the length of the radial segment from $\Sigma$ to radius
$r$, with fixed angular and torus coordinates. Since $A>0$,
$\rho$ is strictly increasing and $\rho(0)=0$. Moreover,
\[
 2a+1=\frac{2(n-4)}{3(n+2)}>0,\qquad
 \rho(r)\sim\frac{\lambda}{2a+1}r^{2a+1}
          \longrightarrow\infty.
\]

We next compare this radial length with the length of an arbitrary
piecewise smooth curve $c:[t_0,t_1]\to M$.
The Euclidean radius $r$ is locally Lipschitz
also at $\Sigma$, so $t\mapsto\rho(r(c(t)))$ is absolutely
continuous. From the radial part $A(r)^2dr^2$ of $g_n$ we obtain,
almost everywhere,
\[
 \left|\frac{d}{dt}\rho(r(c(t)))\right|
 =A(r(c(t)))\left|\frac{d}{dt}r(c(t))\right|
 \leq|\dot c(t)|_{g_n}.
\]
Integrating, for any such curve from $p$ to $q$ we find
\[
 |\rho(r(q))-\rho(r(p))|
 \leq\int_{t_0}^{t_1}
       \left|\frac{d}{dt}\rho(r(c(t)))\right|\,dt
 \leq\ell_{g_n}(c).
\]
Taking the infimum over curves gives
\begin{equation}\label{eq:direct-radial-lipschitz}
 |\rho(r(q))-\rho(r(p))|\leq d_{g_n}(p,q).
\end{equation}
Since $\rho=0$ on
$\Sigma$ and a radial segment realizes this lower bound,
we also obtain
\[
 d_{g_n}(p,\Sigma)=\rho(r(p)).
\]

For each $L\geq0$, let $r_L$ be the unique radius with
$\rho(r_L)=L$. The sublevel set of the radial length is
\[
 \{p:\rho(r(p))\leq L\}
 =\{(X,Y):X^2+Y^2\leq r_L^2\}\times\T^{n-2},
\]
which is compact because both the closed Euclidean disk and the
torus are compact. Now fix a center $o\in M$ and $R>0$.
By \eqref{eq:direct-radial-lipschitz}, every point of the closed
metric ball $\overline B_R^{g_n}(o)$ satisfies
\[
 \rho(r(p))\leq\rho(r(o))+R.
\]
Thus this ball is a closed subset of a compact sublevel set, and
is compact; the Hopf--Rinow theorem ensures that $(M,g_n)$ is complete.

Finally, $S_{g_n}$ is nonconstant, so $g_n$ cannot be Einstein.

\textbf{Step 5.} Finally, we compute the integral of $|S_{g_n}|^p$.
The volume form is
\begin{equation}\label{eq:direct-volume}
 \dv_{g_n}=ABC^{n-2}\,dr\,d\theta\,\dv_{g_0}
 =\lambda^2r(1+r^2)^{a+b+(n-2)c}
 \,dr\,d\theta\,\dv_{g_0}.
\end{equation}
The exponents satisfy
\[
 a+b+(n-2)c=2(2a+1),\qquad
 q_n=2+\frac1{2a+1}.
\]
Consequently, the change of variable $y=1+r^2$ gives
\begin{align*}
 \int_M|S_{g_n}|^p\,\dv_{g_n}
 &=2\pi\lambda^2V_0\int_0^\infty
 r(1+r^2)^{(2-p)(2a+1)}\,dr\\
 &=\pi\lambda^2V_0\int_1^\infty
 y^{(2a+1)(q_n-p)-1}\,dy.
\end{align*}
Since $2a+1>0$, this integral converges precisely when $p>q_n$.
In that case it equals
\[
 \frac{\pi\lambda^2V_0}{(2a+1)(p-q_n)}
 =\frac{16\pi(n-1)V_0}{(n-4)(p-q_n)},
\]
as claimed. Note that, in particular,
\[
 |S_{g_n}|^2\,\dv_{g_n}
 =\lambda^2r\,dr\,d\theta\,\dv_{g_0},
\]
so the energy is infinite.

\end{proof}

\begin{remark}\label{rem:direct-cigar}
When $n=10$, $a=b=-1/3$, $c=1/6$, and $\lambda^2=8$; hence
\[
 g_{10}=(1+r^2)^{1/3}
 \left(8\frac{dr^2+r^2d\theta^2}{1+r^2}+g_0\right).
\]
Note that the two-dimensional metric in the fraction is \emph{Hamilton's cigar}
(see \cite{Hamilton}), and thus the metric is conformal to a gradient steady Ricci soliton.

\end{remark}

\begin{remark}\label{rem:direct-n4}
The same metric formula makes sense at $n=4$, where
$a=-1/2$, $b=-1/6$, $c=1/3$, and $\lambda^2=16/3$.
It gives
\[
 g_4=\frac{16}{3}\frac{dr^2}{1+r^2}
 +\frac{16}{3}r^2(1+r^2)^{-1/3}d\theta^2
 +(1+r^2)^{2/3}g_0
 \quad\text{on }\R^2\times\T^2.
\]
The radial length is
\[
 \rho(r)=\frac4{\sqrt3}\operatorname{arsinh}r\longrightarrow\infty,
\]
so the argument in the previous construction proves completeness. Substitution in
\eqref{eq:direct-Ricci} gives
$\Ric_{g_4}=-g_4/4$ and $S_{g_4}=-1$, i.e. $g_4$ is Einstein and, in particular, it is critical. Its volume form is
$\frac{16}{3}r\,dr\,d\theta\,\dv_{g_0}$, so its volume is infinite.
The non-Einstein examples of
Proposition~\ref{thm:four-complete-examples} require a different
construction.
\end{remark}

\subsection{\texorpdfstring{A six-dimensional example}{A six-dimensional example}}
\label{subsec:six-dimensional-example}

The torus family does not establish the sharpness of the rigidity
endpoint in dimension six, while a  different example does so.

\begin{proposition}\label{prop:n6-sharp-example}
Let $(N^5,\bar h)$ be a closed Einstein manifold with
$\Ric_{\bar h}=-\bar h/10$. Define on $\R\times N$
\begin{equation}\label{eq:n6-g}
 v(t)=\cosh(t/\sqrt{10}),\qquad
 \widehat g_6=v(t)(dt^2+\bar h).
\end{equation}
Then $\widehat g_6$ is smooth, complete, non-Einstein, and critical for
$\cS^2$, with
\begin{equation}\label{eq:n6-scalar}
 S_{\widehat g_6}=-\operatorname{sech}(t/\sqrt{10}).
\end{equation}
Moreover, for every $p>0$,
\begin{equation}\label{eq:n6-integrability}
 S_{\widehat g_6}\in L^p(\R\times N,\widehat g_6)\quad\Longleftrightarrow\quad p>3.
\end{equation}
In particular, the endpoint $P_6=3$ of
Theorem~\ref{thm:rigidity-intro} is optimal.
\end{proposition}

\begin{proof}
Set $h=dt^2+\bar h$, $u=v^{-1}$, and $a=v'/v$.
Since $v''=v/10$ and $\Ric_h=-\bar h/10$, the conformal Ricci
and Hessian formulas for $\widehat g_6=vh$ yield
\[
 \Ric_{\widehat g_6}-u^{-1}\Hess_{\widehat g_6}(u)=-\frac3{20}\,u \widehat g_6,
 \qquad \Delta_{\widehat g_6}u=-\frac1{10}u^2.
\]
Tracing also gives $S_{\widehat g_6}=-u$, as asserted.
As $v\geq1$ and the product metric $h$ is complete, so is $\widehat g_6$.
Its scalar curvature is not constant, hence it is not Einstein.

Put $V_N=\vol_{\bar h}(N)$. Since
$dV_{\widehat g_6}=v^3\,dt\,dV_{\bar h}$, we obtain
\[
 \int_{\R\times N}|S_{\widehat g_6}|^p\,dV_{\widehat g_6}
   =\sqrt{10}\,V_N\int_{\R}\cosh^{3-p}s\,ds.
\]
Using $y=\tanh s$ and then $x=y^2$, the last integral becomes
\begin{equation}\label{eq:n6-exact-integral}
 \int_{-1}^1(1-y^2)^{(p-5)/2}\,dy
   =\int_0^1 x^{-1/2}(1-x)^{(p-5)/2}\,dx.
\end{equation}
At $x=0$, the power $x^{-1/2}$ is integrable; at $x=1$,
integrability is equivalent to $(p-5)/2>-1$, that is, $p>3$.
At $p=3$, the integral over $[-T,T]\times N$ is exactly
$2V_NT$. The distance from $\{0\}\times N$ is
\[
 \rho(t)=\int_0^{|t|}\sqrt{v(\tau)}\,d\tau
      \sim\sqrt{20}\,e^{|t|/(2\sqrt{10})}.
\]
Consequently the integral over the tube $\{\rho\leq R\}$ is
$4\sqrt{10}V_N\log R+O(1)$ as $R\to\infty$.
This specifies the distance with respect to which the endpoint
divergence is logarithmic.
For example, a closed hyperbolic five-manifold, rescaled to sectional curvature
$-1/40$, gives an admissible $N$.
\end{proof}

\begin{remark}
\label{rem:sharp-gradient-examples}
The preceding examples show that the dimensional constant in
\eqref{eq:sharp-intrinsic-gradient-estimate} cannot be decreased
in dimension six or in any dimension $n\geq10$. Indeed, for the torus family of Theorem~\ref{thm:direct-examples},
set $u=-S_{g_n}$; from the radial component of
\eqref{eq:direct-metric} one computes
\[
 H^2=|\nabla u^{-1/2}|^2
   =\frac{(n-4)^2}{24(n-1)(n+2)}
       \frac{r^2}{1+r^2}.
\]
For $n\geq10$, the coefficient equals $K_n$; hence
$\sup_{M_n}H^2=K_n$. In dimension six, the product example
\eqref{eq:n6-g} satisfies
\[
 H^2=|\nabla u^{-1/2}|^2
 =\frac1{40}\tanh^2(t/\sqrt{10}),\qquad
 \sup_{\R\times N}H^2=\frac1{40}=K_6.
\]
All these metrics have $S\in L^p$ for some $p>1$, so they fall
under Proposition~\ref{prop:gradient-estimate}.
In both constructions the sharp constant is approached at an end;
pointwise equality at a finite point is not asserted.
\end{remark}

\section*{Use of artificial intelligence}
The authors used ChatGPT (OpenAI) during the preparation of this manuscript as a supporting tool for mathematical exploration, critical discussion, and editorial revision. The authors take full responsibility for the mathematical content and the final manuscript.

\end{document}